\documentclass[11pt]{amsart}
\usepackage{hyperref}
\usepackage{amscd,amssymb}
\usepackage{amsthm,amsmath,amssymb}
\usepackage{youngtab}
\usepackage{tikz}
\usetikzlibrary{calc}
\usepackage{graphicx}
\usepackage{color}
\usepackage{pdflscape}
\usepackage{bm}
\usepackage{mathrsfs}
\usepackage{verbatim}
\usepackage{adjustbox, collectbox, xkeyval, graphicx, ifpdf, adjcalc}
\newcommand{\ms}{\medskip}

\usepackage{amsmath,amsthm,verbatim, color, amssymb,amsfonts,amscd, graphicx, enumerate}
\usepackage{graphics}
\usepackage{tikz, tikz-cd}
\usetikzlibrary{matrix}
\usepackage[all]{xy}
\usepackage{listings}

\usepackage{comment}

\theoremstyle{plain}
\newtheorem{thm}{Theorem}[section]
\newtheorem{lem}[thm]{Lemma}        
\newtheorem{prop}[thm]{Proposition}
\newtheorem{cor}[thm]{Corollary}

\newtheorem*{thm*}{Theorem}
\newtheorem{abcthm}{Theorem}

\theoremstyle{definition}
\newtheorem{defn}[thm]{Definition} 
\newtheorem*{defn*}{Definition}   
   
\newtheorem*{exam*}{Example}  
\newtheorem{exam}[thm]{Example}  
\newtheorem{rem}[thm]{Remark}

\newcommand{\cd}{\partial}

\newcommand{\lng}{\ell}
\newcommand{\G}{\Gamma}
\newcommand{\Z}{\mathbb{Z}}
\newcommand{\D}{\mathcal{D}}

\newcommand{\Tor}{\operatorname{Tor}}
\newcommand{\id}{\operatorname{id}}
\newcommand{\bt}{\bullet}

\newcommand{\orb}{\mathcal{O}}
\newcommand{\altprod}{\pi}

\definecolor{amethyst}{rgb}{0.6, 0.4, 0.8}
\definecolor{amber}{rgb}{1.0, 0.49, 0.0}
\definecolor{brightmaroon}{rgb}{0.76, 0.13, 0.28}
\definecolor{darkpastelgreen}{rgb}{0.01, 0.75, 0.24}
\usepackage[
paper=a4paper,
headsep=15pt,text={138mm,213mm},centering,includehead
]{geometry}

\usepackage{lipsum}
\address{Institutt for matematiske fag, NTNU Gl\o{}shaugen, Trondheim, Norway}
\email{rachael.boyd@ntnu.no}

\def\subjclassname{\textup{2010} Mathematics Subject Classification}
\expandafter\let\csname subjclassname@1991\endcsname=\subjclassname
\expandafter\let\csname subjclassname@2000\endcsname=\subjclassname

\begin{document}


\title{The low-dimensional homology of finite-rank Coxeter groups}
\author{Rachael Boyd}

\begin{abstract}
	We give formulas for the second and third integral homology of an arbitrary finitely generated Coxeter group, solely in terms of the corresponding Coxeter diagram. The first of these calculations refines a theorem of Howlett, while the second is entirely new and is the first explicit formula for the third homology of an arbitrary Coxeter group.
\end{abstract}

\maketitle

\section{Introduction}

Given a Coxeter group $W$ with finite generating set $S$ and corresponding system~$(W,S)$, denote the associated {Coxeter diagram} by $\D_W$ (see Definitions ~\ref{defn:cox system} and~\ref{defn:cox diagram}).

In this paper, variations on this diagram are defined, and Theorems \ref{thm:INTRO THM A} and \ref{thm:INTRO THM B} below calculate the second and third integral homology for any finite rank Coxeter group~$W$, in terms of zeroth and first cellular homologies of these new diagrams, considered as cell complexes in their own right.

Throughout this paper we will always denote the cyclic group $\Z/n\Z$ as~$\Z_n$. Previously it was known that first and second homology groups of a Coxeter group were isomorphic to $\Z_2^{r_i}$, where $r_i=\operatorname{rank}_{\Z_2}(H_i(W;\Z))$ and both $r_1$ and $r_2$ are known. The computation of $H_1(W;\Z)$ is a straightforward computation of the abelianisation. The computation of $H_2(W;\Z)$ is due to Howlett \cite{Howlett}. A paper of Ihara and Yokonuma \cite{IharaYokonuma} gives results for the second cohomology of certain finite Coxeter groups, with coefficients in $\mathbb{C}^*$. These results agree with Howlett's theorem for the groups in question.

Theorem \ref{thm:INTRO THM A} below gives a refinement of Howlett's theorem by introducing a naturality statement. The method of proof is new and uses a spectral sequence argument. Theorem \ref{thm:INTRO THM B} is the first explicit formula for $H_3(W;\Z)$ and extends the same method. This method could be extended to produce computations of higher homologies, the drawback being that the differentials in the spectral sequence become more difficult to handle as the homological degree increases. Terms that we use while stating our results below will be defined in Section \ref{sec:coxeter def and egs}.

\subsection{Second homology}

Given a diagram $\D$, let $E(\D)$ and $V(\D)$ be the set of edges and set of vertices of $\D$ respectively. Let~$\D_W$ be the Coxeter diagram corresponding to Coxeter system~$(W,S)$. Then~$V(\D_W)=S$ and to every pair~$s\neq t \in S$ there is an associated label~$m(s,t) \in \mathbb{N}\cup \infty$.

\begin{defn}\label{def:diagrams for h2}
	We introduce three new diagrams: $\D_{\operatorname{odd}}$, $\D_{\operatorname{even}}$ and $\D_{\bt \bt}$. 
\begin{itemize}
	\item Let $\D_{\operatorname{odd}}$ be the diagram with $V(\D_{\operatorname{odd}})=S$ and
	$$
	e(s,t) \in E(\D_{\operatorname{odd}}) \iff m(s,t) \text{ is odd}.
	$$
	\item Let $\D_{\operatorname{even}}$ be the diagram with  $V(\D_{\operatorname{even}})=S$ and
	$$
	e(s,t) \in E(\D_{\operatorname{even}}) \iff 2\neq m(s,t) \text{ is even}.
	$$ 
	\item Let $\D_{\bullet \bullet }$ be the diagram with $$V(\D_{\bullet \bullet })=\{\{s,t\} \mid s, t \in \, S, \, m(s,t)=2\}$$ $$  e(\{s_1,t_1\},\{s_2,t_2\}) \in E(\D_{\bullet \bullet }) \iff s_1=s_2\text{ and }m(t_1,t_2) \text{ is odd.} $$
\end{itemize}
\end{defn}

\begin{abcthm} \label{thm:INTRO THM A}
	Given a finite rank Coxeter system $(W,S)$, there is a natural isomorphism
	$$
	H_2(W;\Z) \overset{}{\cong}H_0(\D_{\bullet \bullet}; \Z_2) \oplus \Z_2[E(\D_{\operatorname{even}})] \oplus H_1(\D_{\operatorname{odd}};\Z_2)
	$$
	where in the first and final term of the right-hand side the diagrams are considered as 1-dimensional cell complexes.
\end{abcthm}

\begin{rem}
 Computing the rank of the right-hand side recovers Howlett's theorem \cite{Howlett}.
\end{rem}

Consider the category where the objects are Coxeter systems and the morphisms are full inclusions (Definition \ref{def:full inclusion}), then group homology acts as a functor to the category of abelian groups. The right-hand side of the isomorphism in Theorem~\ref{thm:INTRO THM A} assigns to a Coxeter diagram $\D_W$ the three new diagrams $\D_{\operatorname{odd}}$, $\D_{\operatorname{even}}$ and $\D_{\bt \bt}$ and furthermore assigns to these diagrams an abelian group. The total outcome is again a functor to abelian groups. Naturality says that the isomorphism of the statement is a natural isomorphism of functors.

\subsection{Third homology}

To state this theorem we introduce four new diagrams.

\begin{defn}\label{def:diagrams for h3} Let~$\D_W$ be a Coxeter diagram corresponding to the Coxeter system~$(W,S)$.
\begin{itemize}
	\item Let $\D_{{\bf A}_2}$ be the diagram with $$V(\D_{{\bf A}_2})=\{\{s,t\} \mid s, t \in \, S, \, m(s,t)=3\}$$  $$e(\{s_1,t_1\},\{s_2,t_2\}) \in E(\D_{{\bf A}_2}) \iff s_1=s_2\text{ and } m(t_1,t_2)=2.$$
	\item Let $\D{\begin{tikzpicture}[scale=0.08, baseline={(0,0)}]
		\draw[fill= black] (12.5,0) circle (0.6) ;
		\draw (17.5,1.4) node {{\tiny 2r}};
		\draw[line width=1] (15,0) -- (20,0);
		\draw[fill= black] (15,0) circle (0.6);
		\draw[fill= black] (20,0) circle (0.6); 
		\end{tikzpicture}}$ be the diagram with $$V(\D{\begin{tikzpicture}[scale=0.08, baseline={(0,0)}]
		\draw[fill= black] (12.5,0) circle (0.6) ;
		\draw (17.5,1.4) node {{\tiny 2r}};
		\draw[line width=1] (15,0) -- (20,0);
		\draw[fill= black] (15,0) circle (0.6);
		\draw[fill= black] (20,0) circle (0.6); 
		\end{tikzpicture}})=\{\{s,t, u\} \mid s, t, u \, \in \, S, \, m(s,t)=m(s,u)=2, m(t,u)=2r \text{ is even}\}$$	
	 $$e(\{s_1,t_1, u_1\},\{s_2,t_2, u_2\}) \in E(\D{\begin{tikzpicture}[scale=0.08, baseline={(0,0)}]
	 	\draw[fill= black] (12.5,0) circle (0.6) ;
	 	\draw (17.5,1.4) node {{\tiny 2r}};
	 	\draw[line width=1] (15,0) -- (20,0);
	 	\draw[fill= black] (15,0) circle (0.6);
	 	\draw[fill= black] (20,0) circle (0.6); 
	 	\end{tikzpicture}}) \iff \begin{array}{l}
 	t_1=t_2, \, u_1=u_2 \\m(s_1,s_2)\text{ is odd.}
 \end{array}$$ 
	\item Let $\D_{\bf{A}_3}$ be the diagram with $$ V(\D_{\bf{A}_3})=\{\{s,t,u\} \mid s, t,u \in \, S, \,  m(s,t)=m(t,u)=3 \text{ and }m(s,u)=2 \}$$
	$$ e(\{s_1,t_1, u_1\},\{s_2,t_2,u_2\}) \in E(\D_{\bf{A}_3}) \iff t_1=t_2,\, u_1=u_2, \, m(s_1,s_2)=2.$$ 	
	\item Let $\D^\square_{\bt \bt}$ be the CW-complex formed from the diagram $\D_{\bullet \bullet }$ by attaching a $2$-cell to every square. 
\end{itemize}
\end{defn}

\begin{abcthm}\label{thm:INTRO THM B}
	Given a finite rank Coxeter system $(W,S)$ there is an isomorphism
\begin{eqnarray*}
	H_3(W; \Z) &\cong& H_0(\D_{\operatorname{odd}};\Z_2)  \oplus H_0(\D_{{\bf A}_2};\Z_3) \oplus \big(\underset{3<m(s,t)< \infty}{\oplus}\Z_{m(s,t)} \big)\\&&\oplus H_0(\D{\begin{tikzpicture}[scale=0.08, baseline={(0,0)}]
		\draw[fill= black] (12.5,0) circle (0.6) ;
		\draw (17.5,1.4) node {{\tiny 2r}};
		\draw[line width=1] (15,0) -- (20,0);
		\draw[fill= black] (15,0) circle (0.6);
		\draw[fill= black] (20,0) circle (0.6); 
		\end{tikzpicture}};\Z_2) \oplus \big(\underset{\substack{{\scriptscriptstyle W({\bf H}_3)\subseteq W}\\{\scriptscriptstyle W({\bf B}_3)\subseteq W}}}{\oplus}\Z_2 \big)\\&&\oplus \big(H_0(\D_{\bf{A}_3};\Z_2) \bigcirc H_0(\D_{\bullet \bullet };\Z_2)\big)  	\oplus H_1(\D^\square_{\bt \bt};\Z_2)
\end{eqnarray*}
where each diagram is viewed as a cell complex.
In this equation, $\bigcirc$ denotes a known non-trivial extension of $H_0(\D_{\bf{A}_3};\Z_2)$ by $H_0(\D_{\bullet \bullet };\Z_2)$ fully described via an extension matrix $X_W$ from Definition \ref{def:extension matrix}.		
\end{abcthm}

We note that the unpublished PhD thesis of Harris~\cite{Harristhesis} contains an independent computation of the third integral homology of a Coxeter group, which differs from Theorem~\ref{thm:INTRO THM B} in many cases.

The finite Coxeter groups were classified in the 1930s by Coxeter \cite{Coxeter}. This classification is described in Theorem \ref{prop:classification of finite Coxeter}. We use Theorems \ref{thm:INTRO THM A} and \ref{thm:INTRO THM B} to calculate the second and third integral homology of the finite Coxeter groups, and give the results in Appendix \ref{appendix:table}. 

\subsection{Outline of proof}
Given a Coxeter system $(W,S)$ these results arise from the computation of the {isotropy spectral sequence} for a contractible $CW$-complex upon which the Coxeter group $W$ acts, called the {Davis complex}. Cells in the Davis complex correspond to finite Coxeter groups that appear in $W$, the {spherical subgroups}. These have Coxeter systems $(W_T,T)$ where $W_T$ is a finite Coxeter group and $T\subseteq S$. The set of $T\subseteq S$ which generate spherical subgroups of a fixed Coxeter group $W$ is denoted $\mathcal{S}$.

The isotropy spectral sequence abuts to the integral homology of $W,$ and the $E^1$ terms are given by the sums of twisted homologies of the spherical subgroups $W_T$ of $W$, for $T$ a given size. 

$$
E^1_{p,q} =\bigoplus_{\substack{T \in \mathcal{S}\\ | T | = p}} H_q(W_T; \Z_T) \Rightarrow H_{p+q}(W;\Z).
$$

For the proof of Theorem \ref{thm:INTRO THM A} the groups on the $E^1$ terms and $d^1$ differential of the spectral sequence are simple to compute. We see there are no further differentials that will affect the diagonal corresponding to $H_2(W;\Z)$ on the $E^\infty$ page, so the limiting terms are equal to the $E^2$ terms. There is only one non-zero term on the diagonal so there are no possible extension problems and Theorem \ref{thm:INTRO THM A} follows.

For Theorem \ref{thm:INTRO THM B}, the computation of the ~$E^1$ terms relies heavily on a free resolution for Coxeter groups, described by De Concini and Salvetti~\cite{DCSal}. The computer algebra package PyCox~\cite{Geck} is used (though not strictly necessary) to complete some of the longer calculations required.

In order to apply the $d^1$ differential to computations using this resolution, a chain map between resolutions is computed in the required degrees. Using these tools, the $E^2$ page of the spectral sequence on the diagonal corresponding to $H_3(W;\Z)$ is computed. Following this, we use a variety of techniques to prove that all further differentials to and from this diagonal are in fact zero. This includes defining a pairing for the isotropy spectral sequence.

The possible extension problems arising on the limiting page at this diagonal are treated by considering representing subgroups of $W$ for each class and mapping between the corresponding spectral sequences. From these computations we note there is only one non-trivial extension and thus Theorem \ref{thm:INTRO THM B} follows.

\subsection{Organisation of paper}
We start with background on Coxeter groups and an introduction to the Davis complex $\Sigma_W$ of $W$ in Section \ref{sec:coxeter def and egs}. We then introduce the isotropy spectral sequence in Section \ref{sec:isotropy spectral seq}, and prove some associated desired results. Following this, Section \ref{sec:h2 calc coxeter} proves Theorem \ref{thm:INTRO THM A} and Section \ref{SEC:H3 FOR COXETER} proves Theorem \ref{thm:INTRO THM B}.

\subsection{Acknowledgements}

I would like to thank my PhD advisor Richard Hepworth for his wonderful teaching, our many maths conversations and his unconditional support. I would also like to thank my thesis examiners Mark Grant and Ulrike Tillmann for their careful reading and related comments and suggestions, and Markus Szymik for helpful feedback. Finally I would like to thank the anonymous referee for valuable comments.

\section{Coxeter groups}\label{sec:coxeter def and egs}

\noindent This expository section follows~\cite{Davis}.

\begin{defn} \label{defn:cox system}
	A \emph{Coxeter matrix} on a finite set $S$ is a symmetric~$S \times S$ matrix~$M$ with entries~$m(s,t)$ in $\mathbb{N} \cup \{\infty\}$ for $s, t$ in $S$. This matrix must satisfy~$m(s,t)=1$ if and only if~$s=t$ and $m(s,t)=m(t,s)$ must be greater than 1 when $s\neq t$.
	A Coxeter matrix $M$ has an associated \emph{Coxeter group} $W$, with presentation
	$$
	W=\langle S \, \vert \,  (st)^{m(s,t)}=e\rangle.
	$$
	\noindent We call $(W,S)$ a \emph{Coxeter system}, and we call $|S|$ the \emph{rank} of $(W,S)$. We adopt the convention that $(W,\emptyset)$ is the trivial group.
\end{defn}

\begin{rem}
	Note that the condition $m(s,s)=1$ implies that all generators of the group are involutions i.e.~$s^2=e$ for all $s$ in $S$.
\end{rem}

\begin{defn}
	Define the \emph{length function }on a Coxeter system~$(W,S)$ to be the function~$\lng: W \to \mathbb{N}$ which maps~$w$ in~$W$ to the minimal word length required to express~$w$ in terms of the generators in~$S$. That is, we set $\lng(e)=0$, and if $w\neq e$ then there exists a minimal~$k\geq 1$ such that~$w=s_1\cdots s_k$ for~$s_i$ in~$S$.
\end{defn}

\begin{defn}\label{def: pi(s,t;m)}
	For $k \in \mathbb{N}$, define $\altprod(a,b ; k)$ to be the word of length $k$, given by the alternating product of $a$ and $b$ i.e. \vspace{-.35cm}
	$$
	\altprod(a,b ; k)=\overbrace{abab\ldots }^\text{length k}.
	$$
\end{defn}

\begin{rem}\label{rem:second pres for W}
	When $m(s,t)\neq \infty$, the relations $(st)^{m(s,t)}=e$ can be rewritten as
	$$
	\altprod(s,t;m(s,t))=\altprod(t,s;m(s,t)).
	$$
\end{rem} 

\begin{defn}\label{defn:cox diagram}
Given a Coxeter system $(W,S)$, the associated \emph{Coxeter diagram}, denoted $\D_W$, is a labelled graph with vertices indexed by the generating set $S$. Edges are drawn between the vertices corresponding to $s$ and $t$ in $S$ when $m(s,t)\geq~3$ and labelled with $m(s,t)$ when $m(s,t)\geq 4$ (or $\infty$). When the diagram $\D_W$ is connected, $W$ is called an \emph{irreducible} Coxeter system.
\end{defn}

\begin{thm}[{Classification of finite Coxeter groups, Coxeter \cite{Coxeter}}]\label{prop:classification of finite Coxeter}
	A Coxeter system is finite (i.e.~gives rise to a finite Coxeter group) if and only if it is the (direct) product of finitely many finite irreducible Coxeter systems. The following is a complete list of the diagrams corresponding to finite irreducible Coxeter systems, and therefore classifies finite Coxeter groups.
	
	\ms
	
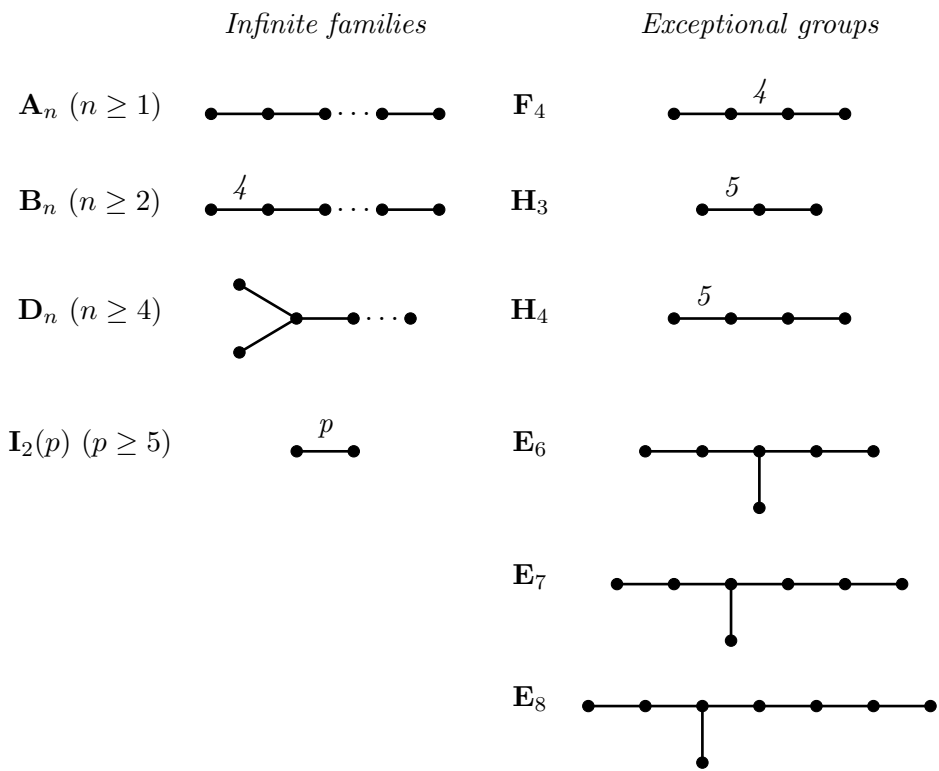
\centerline{	\xymatrix@R=3mm@C=2mm {	
			& \textrm{Infinite families} & & & \textrm{Exceptional groups}\\	
			{\bf A}_n \,\, (n\geq 1) & \begin{tikzpicture}[scale=0.15, baseline=0]
			\draw[fill= black] (5,0) circle (0.5);
			\draw[line width=1] (5,0) -- (15,0);
			\draw[fill= black] (10,0) circle (0.5);
			\draw[fill= black] (15,0) circle (0.5);		
			\draw (17.5,0) node {$\ldots$};
			\draw[fill= black] (20,0) circle (0.5);
			\draw[line width=1] (20,0) -- (25,0);				
			\draw[fill= black] (25,0) circle (0.5);
			\end{tikzpicture} 
			& & {\bf F}_4 & \begin{tikzpicture}[scale=0.15, baseline=0]
			\draw[fill= black] (5,0) circle (0.5);
			\draw[line width=1] (5,0) -- (20,0);
			\draw[fill= black] (10,0) circle (0.5);
			\draw[fill= black] (15,0) circle (0.5);
			\draw[fill= black] (20,0) circle (0.5);
			\draw (12.5,2) node {4};
			\end{tikzpicture}  \\
			{\bf B}_n \,\,(n\geq2) & \begin{tikzpicture}[scale=0.15, baseline=0]
			\draw[fill= black] (5,0) circle (0.5);
			\draw[line width=1] (5,0) -- (15,0);
			\draw[fill= black] (10,0) circle (0.5);
			\draw[fill= black] (15,0) circle (0.5);		
			\draw (17.5,0) node {$\ldots$};
			\draw[fill= black] (20,0) circle (0.5);
			\draw[line width=1] (20,0) -- (25,0);				
			\draw[fill= black] (25,0) circle (0.5);
			\draw (7.5,2) node {4};
			\end{tikzpicture} 
			& & {\bf H}_3 & \begin{tikzpicture}[scale=0.15, baseline=0]
			\draw[fill= black] (5,0) circle (0.5);
			\draw[line width=1] (5,0) -- (15,0);
			\draw[fill= black] (10,0) circle (0.5);
			\draw[fill= black] (15,0) circle (0.5);
			\draw (7.5,2) node {5};
			\end{tikzpicture}  \\
			{\bf D}_n \,\, (n\geq 4) & \begin{tikzpicture}[scale=0.15, baseline=0]
			\draw[fill= black] (5,-3) circle (0.5);
			\draw[line width=1] (5,-3) -- (10,0);
			\draw[line width=1] (5,3) -- (10,0);
			\draw[line width=1] (15,0) -- (10,0);
			\draw[fill= black] (5,3) circle (0.5);
			\draw[fill= black] (10,0) circle (0.5);
			\draw[fill= black] (15,0) circle (0.5);
			\draw (17.5,0) node {$\ldots$};				
			\draw[fill= black] (20,0) circle (0.5);
			\end{tikzpicture}  
			& & {\bf H}_4 & \begin{tikzpicture}[scale=0.15, baseline=0]
			\draw[fill= black] (5,0) circle (0.5);
			\draw[line width=1] (5,0) -- (20,0);
			\draw[fill= black] (10,0) circle (0.5);
			\draw[fill= black] (15,0) circle (0.5);
			\draw[fill= black] (20,0) circle (0.5);
			\draw (7.5,2) node {5};
			\end{tikzpicture}  \\
			{\bf I}_2(p) \,\, (p\geq 5)& \begin{tikzpicture}[scale=0.15, baseline=0]
			\draw[fill= black] (5,0) circle (0.5);
			\draw (7.5,2) node {p};
			\draw[line width=1] (5,0) -- (10,0);			
			\draw[fill= black] (10,0) circle (0.5);
			\end{tikzpicture} 
			& &{\bf E}_6 & \begin{tikzpicture}[scale=0.15, baseline=0]
			\draw[fill= black] (5,0) circle (0.5);
			\draw[line width=1] (5,0) -- (25,0);
			\draw[fill= black] (10,0) circle (0.5);
			\draw[fill= black] (15,0) circle (0.5);
			\draw[fill= black] (20,0) circle (0.5);		
			\draw[fill= black] (25,0) circle (0.5);			
			\draw[fill= black] (15,-5) circle (0.5);
			\draw[line width=1] (15,0) -- (15,-5);
			\end{tikzpicture}  \\ 
			& & & {\bf E}_7 & \begin{tikzpicture}[scale=0.15, baseline=0]
			\draw[fill= black] (5,0) circle (0.5);
			\draw[line width=1] (5,0) -- (30,0);
			\draw[fill= black] (10,0) circle (0.5);
			\draw[fill= black] (15,0) circle (0.5);
			\draw[fill= black] (20,0) circle (0.5);		
			\draw[fill= black] (25,0) circle (0.5);		
			\draw[fill= black] (30,0) circle (0.5);		
			\draw[fill= black] (15,-5) circle (0.5);
			\draw[line width=1] (15,0) -- (15,-5);
			\end{tikzpicture}   \\
			& & & {\bf E}_8 & \begin{tikzpicture}[scale=0.15, baseline=0]
			\draw[fill= black] (5,0) circle (0.5);
			\draw[line width=1] (5,0) -- (35,0);
			\draw[fill= black] (10,0) circle (0.5);
			\draw[fill= black] (15,0) circle (0.5);
			\draw[fill= black] (20,0) circle (0.5);		
			\draw[fill= black] (25,0) circle (0.5);		
			\draw[fill= black] (30,0) circle (0.5);		
			\draw[fill= black] (35,0) circle (0.5);		
			\draw[fill= black] (15,-5) circle (0.5);
			\draw[line width=1] (15,0) -- (15,-5);
			\end{tikzpicture}  \\ 
	}}
\end{thm}

\paragraph{\bf{Notation.}} Throughout this paper, for ease of notation we may write~${\bf I}_2(2)$,~${\bf I}_2(3)$ and ${\bf I}_2(4)$ instead of~${\bf A}_1\times {\bf A}_1$, ${\bf A}_2$ and~${\bf B}_2$ respectively. Whenever we write~${\bf I}_2(p)$ we will specify for which~$p$ the result corresponds.

\begin{defn}
	We say that a finite irreducible Coxeter group $W$ is \emph{of type} $\D$ if its corresponding diagram is given by $\D$, and we denote this Coxeter group $W(\D)$.
\end{defn}

\begin{rem}
	The Coxeter group of type ${\bf A}_n$, or $W({\bf A}_n)$, is isomorphic to the symmetric group $S_{n+1}$ and the Coxeter group of type ${\bf I}_2(p)$, or $W({\bf I}_2(p))$, is isomorphic to the dihedral group $D_{2p}$. Similarly, the Coxeter group of type ${\bf B}_n$, or $W({\bf B}_n)$, is isomorphic to the signed permutation group $\Z_2 \wr S_{n}$ and $W({\bf D}_n)$, is isomorphic to an index two subgroup of $W({\bf B}_n)$, such that the signs in each permutation multiply to $+1$.
\end{rem}

\subsection{Products and subgroups}\label{sec:cox products and subgroups}

Consider two Coxeter systems $(U, S_U)$ and $(V, S_V)$ and denote by $\D_U \sqcup \D_V$ the diagram created by placing $\D_U$ and $\D_V$ beside each other, disjointly. 

\begin{lem}
	The diagram $\D_U \sqcup \D_V$ defines a Coxeter group $W\cong U\times V$, with diagram $\D_W=\D_U \sqcup D_V$ and generating set $S_W:=S_U \cup S_V$.
\end{lem}

\begin{defn}\label{def:full inclusion}
	A map $\iota\colon\D_U {\rightarrow} \D_{W}$ of Coxeter diagrams is a \emph{full inclusion} if~$\iota\colon~U\to W$ is injective and for every $s, t \in U$, $m(\iota(s),\iota(t))=m(s,t)$. In this setting we call $\D_U$ a \emph{full subdiagram} of $\D_W$.
\end{defn} 

\begin{defn}
	Let $(W,S)$ be a Coxeter system. For each $T\subseteq S$ denote by~$W_T$ the subgroup of $W$ generated by $T$. We call subgroups that arise in this way \emph{parabolic subgroups}.
\end{defn}

\begin{prop}[{\cite[4.1.6.(i)]{Davis}}] \label{prop:fullinclusions}
	For $W_T$ a parabolic subgroup, $(W_T, T)$ is a Coxeter system in its own right, and defines a full inclusion $\D_{W_T} \hookrightarrow \D_W$. Similarly, a full inclusion corresponds to a parabolic subgroup.
\end{prop}

The next result concerns cosets of parabolic subgroups. Let $(W,S)$ be a Coxeter system, and $T, T'$ be subsets of $S$.

\begin{lem} [{\cite[4.3.1]{Davis}}] \label{lem:shortest coset element coxeter}
	There is a unique element of minimal length in the double coset $W_T w W_{T'}$.
\end{lem}

\begin{defn}[{\cite[4.3.2]{Davis}}]\label{def:(-,-) reduced coxeter}
	We say an element $w$ in $W$ is \emph{$(T,T')$-reduced} if $w$ is the shortest element in $W_T w W_{T'}$.
\end{defn}

\begin{rem}	
	Given the parabolic subgroup $W_T$ in $W$, $w$ in $W$ is $(T, \emptyset)$-reduced if $\lng(tw)=\lng(t)+\lng(w)=1+\lng(w)$ for all $t$ in $T$. Note that this implies~$w$ cannot be written in such a way that it starts with any letter in $T$. Likewise we say $w$ in $W$ is $( \emptyset, T)$-reduced if $\lng(wt)=\lng(w)+1$ for all $t$ in $T$.
\end{rem}

\begin{defn}
	A finite parabolic subgroup is called a \emph{spherical subgroup}.
\end{defn} 

Since the diagrams of parabolic subgroups appear as full subdiagrams of the Coxeter diagram, given a Coxeter system $(W,S)$ we identify its spherical subgroups via occurrences of the irreducible diagrams from Theorem \ref{prop:classification of finite Coxeter} in $\D_W$, and disjoint unions of such diagrams.
\begin{defn}\label{def:mathcalS}
	Given a Coxeter stem $(W,S)$, we denote by $\mathcal{S}$ the set of all subsets of $S$ which generate spherical subgroups of $W$, i.e.
	$$
	\mathcal{S}=\{T \subseteq S \mid W_T \text{ is finite}\}.
	$$
\end{defn}

\subsection{The Davis complex}\label{sec:Davis complex}

In this section we introduce the Davis complex for a Coxeter group.

\begin{defn}
	A coset of a spherical subgroup is called a \emph{spherical coset}. For a Coxeter system $(W,S)$ and a subgroup $W_T$ we denote the set of cosets as follows:
	$$
	W/W_T=\{wW_T \,\vert \, w \in W\}.
	$$
	The \emph{poset of spherical cosets} is denoted $W\mathcal{S}$:
	$$
	W\mathcal{S}=\bigcup_{T \in \mathcal{S}} \{W/W_T\}.
	$$
	where $W\mathcal{S}$ is partially ordered by inclusion. The group $W$ acts on the poset $W\mathcal{S}$ by left multiplication and the quotient poset is $\mathcal{S}$.
\end{defn} 

\begin{lem}[{\cite[4.1.6.(iii)]{Davis}}]\label{lem:coset inclusion coxeter}
	Given $T$ and $U$ in $\mathcal{S}$ and $w$ and $v$ in $W$, the cosets $wW_U$ and $vW_T$ satisfy $wW_U\subseteq vW_T$ if and only if $w^{-1}v \in W_T$ and $U \subseteq T$.
\end{lem}

\begin{defn}[{\cite[7.2]{Davis}}]
	One can associate to a Coxeter system $(W,S)$ a CW-complex called the \emph{Davis complex}. This is denoted $\Sigma_W$ and is the geometric realisation of the poset $W\mathcal{S}$. That is every spherical coset $wW_T$ is realised as a vertex or $0$-cell, and for every ordered chain of $(p+1)$ spherical cosets there is a~$p$-cell in the Davis complex:
	$$
	w_0W_{T_0}\subset w_1W_{T_1}\subset w_2W_{T_2}\subset \cdots \subset w_pW_{T_p}
	$$
	where $w_i$ is in $W$ and $T_i$ is in $\mathcal{S}$ for all $0\leq i \leq p$. The associated Coxeter group $W$ acts on the Davis complex by left multiplication on the cosets.
\end{defn}

\begin{defn}
	For every finite Coxeter group $W$ with generating set $S$, one can define a \emph{canonical representation} of the Coxeter group $W$  on $\mathbb{R}^n$, where $n=|S|$ (see \cite[Section~6.12]{Davis} for details). Given this representation we define the \emph{Coxeter polytope}, or \emph{Coxeter cell} of $W$ to be the convex hull of the orbit of a generic point~$x$ in $\mathbb{R}^n$ under the $W$-action. This polytope has dimension $n=|S|$, and we denote it~$C_W$. A detailed definition can be found in \cite[Section 7.3]{Davis}.
\end{defn}

\begin{prop}[{\cite[Lemma 7.3.3]{Davis}}]
	If $W$ is a finite Coxeter group then~$\Sigma_W$ is homeomorphic to the barycentric subdivision of the Coxeter cell $C_W$. 
\end{prop}

\begin{defn}	
	A coarser cell structure can be given to~$\Sigma_W$:~consider only those spherical cosets which are present as subsets of a chosen coset~$wW_T$ and denote this~$W\mathcal{S}_{\leq wW_T}$. The realisation of~$W\mathcal{S}_{\leq wW_T}$ is a subcomplex of~$\Sigma_W$. In fact~$|W\mathcal{S}_{\leq wW_T}| \cong |W_T\mathcal{S}_T|$ where~$\mathcal{S}_T$ denotes the set of spherical subsets of $T$. Since $W_T$ is finite, the realisation of~$W_T\mathcal{S}_T$ is homeomorphic to the barycentric subdivision of its Coxeter cell~$C_{W_T}$. Therefore the realisation is homeomorphic to a disk, i.e.~$|W_T\mathcal{S}_T|\cong \mathbb{D}^{|T|}$. The cell structure on~$\Sigma_W$ is therefore given by associating to the subcomplex~$W\mathcal{S}_{\leq wW_T}$ its corresponding Coxeter cell: a~$p$-cell where~$p=|T|$. The~$0$-cells are given by cosets of the form~$W\mathcal{S}_{\leq wW_{\emptyset}}$, i.e.~the set~$\{wW_{\emptyset}| w \in W\}$, and therefore associated to elements of~$W$ (recall $W_{\emptyset}=\{e\}$). By Lemma \ref{lem:coset inclusion coxeter} a set of vertices~$X$ will define a~$p$-cell precisely when~$X= wW_T$ for~$T \in \mathcal{S}$ and $|T|=p$. There is an action of~$W$ on the cells of~$\Sigma_W$ given by left multiplication, and this makes~$\Sigma_W$ into a~$W$-complex in the sense of~\cite{Brown}. The stabiliser of a~$p$-cell~$wW_T$ under this action is the finite subgroup~$wW_Tw^{-1}$ and upon identification of the cell~$wW_T$ with~$C_{W_T}$ this acts by reflections in the usual way.
\end{defn}

\noindent We use the following results concerning the Davis complex in this paper.

\begin{prop}[{\cite[8.2.13]{Davis}}]\label{prop: Davis is contractible}
	For any Coxeter group $W$, $\Sigma_W$ is contractible.
\end{prop}

\begin{lem}[{\cite[7.4.4]{Davis}}] \label{lem:product of Davis}
	Suppose $W$ and $S$ decompose as~$W=U\times V$ and~$S=S_U\cup S_V$. Then~$\mathcal{S}=\mathcal{S}_U \times \mathcal{S}_W$ and~$\Sigma_W=\Sigma_U \times \Sigma_V$ is an isomorphism of CW-complexes, provided we use the coarser cell structure.
\end{lem}

\section{The isotropy spectral sequence}\label{sec:isotropy spectral seq}

We give explicit formulas for the terms on the $E^1$ page of the isotropy spectral sequence for the Davis complex, as well as the $d^1$ differential, which is induced by a transfer map. We also introduce a pairing for the isotropy spectral sequence of the Davis complex in Section \ref{sec:pairings}.

\subsection{Isotropy spectral sequence for the Davis Complex}\label{sec:isotropy spectral seq for davis}

We consider the isotropy spectral sequence for a Coxeter system $(W,S)$ and related~$\Sigma_W$. Recall the definition of the isotropy spectral sequence from \cite{Brown}[VII, Equation(7.10)]. For more background see \cite{Boydthesis} or \cite{Brown}.

Consider the action of~$W$ on~$\Sigma_W$ and denote the stabiliser of a cell~$\sigma$ by $W_\sigma$. Denote the orientation module of~$\sigma$ by~$\Z_\sigma$.  We consider the isotropy spectral sequence for integral homology. 

\begin{lem}\label{lem:isotropy for davis}
	Under the $W$-action on $\Sigma_W$, a set of orbit representatives of $p$-cells is
	$$
	\orb_p=\{eW_T \mid T \in \mathcal{S} \, ,\,  |T|=p\}.
	$$
	The stabilizer of a cell $\sigma=eW_T$ is $W_\sigma=W_T$ and the action of an element $w$ of~$W_T$ on~$\Z_\sigma$ is the identity if $\lng(w)$ is even, or negation if $\lng(w)$ is odd.
	\begin{proof}
		Recall that each~$p$-cell of~$\Sigma_W$ is represented by a spherical coset~$wW_T$ where~$|T|=p$, and the vertices of the cell are given by the set $\{vW_\emptyset| v \in wW_T\}$. The group $W$ acts by left multiplication and so we can choose the orbit representatives of $p$-cells to be the cosets $eW_T=W_T$ where~$|T|=p$ and~$T$ is in~$\mathcal{S}$.
		The stabiliser of a cell represented by $W_T$ is~$W_T$ itself. Every element in the generating set $T$ of~$W_T$ acts on the cell by reflection, reversing the orientation of the cell. The action of an element of~$W_T$ on the orientation module will therefore be the identity if the element has even length, or negation if the element has odd length. 
	\end{proof}
\end{lem}

Recall that the Davis complex is contractible (Proposition \ref{prop: Davis is contractible}) and hence acyclic. Then under the choices of Lemma \ref{lem:isotropy for davis}, the isotropy spectral sequence~is

$$
E^1_{p,q}=H_q(W; C_{p}(\Sigma_W, \Z))=\bigoplus_{\sigma \in \orb_p} H_q(W_\sigma; \Z_\sigma) =\bigoplus_{\substack{T \in \mathcal{S}\\ | T | = p}} H_q(W_T; \Z_T) \Rightarrow H_{p+q}(W;\Z)
$$

\noindent since $\Z_\sigma\otimes \Z\cong \Z_\sigma$, which we write as $\Z_T$ for the orientation module of the cell~$W_T$. This gives $E^1$ page as shown in Figure \ref{fig:isotropy ss}.
	\begin{figure}
		\begin{tikzpicture}
		\matrix (m) [matrix of math nodes,
		nodes in empty cells,nodes={minimum width=5ex,
			minimum height=5ex,outer sep=-5pt},row sep=1ex]{
			3 	& 	\color{teal}{H_3(W_\emptyset; \Z_\emptyset)}	& \xleftarrow{d^1} \underset{t \in S}{\oplus} H_3(W_t; \Z_t)  &  \xleftarrow{d^1} \underset{| T | =2}{\underset{T \in \mathcal{S}}{\oplus}}  H_3(W_T; \Z_T)&  \xleftarrow{d^1} \underset{| T | =3}{\underset{T \in \mathcal{S}}{\oplus}}  H_3(W_T; \Z_T) \\
			2 	& 	\color{red}{H_2(W_\emptyset; \Z_\emptyset)}	& \xleftarrow{d^1} \color{teal}{\underset{t \in S}{\oplus} H_2(W_t; \Z_t)}  &  \xleftarrow{d^1} \underset{| T | =2}{\underset{T \in \mathcal{S}}{\oplus}}  H_2(W_T; \Z_T)&  \xleftarrow{d^1} \underset{| T | =3}{\underset{T \in \mathcal{S}}{\oplus}}  H_2(W_T; \Z_T)  \\
			1   &  H_1(W_\emptyset; \Z_\emptyset)  & \xleftarrow{d^1}\color{red}{\underset{t \in S}{\oplus} H_1(W_t; \Z_t)}  & \xleftarrow{d^1} \color{teal}{\underset{| T | =2}{\underset{T \in \mathcal{S}}{\oplus}}  H_1(W_T; \Z_T)}&  \xleftarrow{d^1} \underset{| T | =3}{\underset{T \in \mathcal{S}}{\oplus}}  H_1(W_T; \Z_T) \\
			0   &  	H_0(W_\emptyset; \Z_\emptyset) 	& \xleftarrow{d^1}\underset{t \in S}{\oplus} H_0(W_t; \Z_t) & \xleftarrow{d^1} \color{red}{\underset{| T | =2}{\underset{T \in \mathcal{S}}{\oplus}}   H_0(W_T; \Z_T) }  & \xleftarrow{d^1} \color{teal}{\underset{| T | =3}{\underset{T \in \mathcal{S}}{\oplus}}  H_0(W_T; \Z_T)} \\
			\quad\strut &   0  			&  1  &  2  &  3 \strut \\};
		\draw[thick] (m-1-1.east) -- (m-5-1.east) ;
		\draw[thick] (m-5-1.north) -- (m-5-5.north) ;
		\end{tikzpicture}
		\caption{The $E^1$ page of the isotropy spectral sequence for the Davis complex}
		\label{fig:isotropy ss}
	\end{figure}
\noindent The zeroth column only has one summand, since only the empty set satisfies the criteria of generating a spherical subgroup and having size zero. For the first column, note that all generators in $S$ generate a cyclic group of order two. Denote the subgroup generated by~$s$ in~$S$ by~$W_s$.

We denote the $d^1$ differential component restricted to the $H_q(W_T; \Z_T)$ component in the source and projected to the $H_q(W_U; \Z_U)$ component in the target by $d^1_{T,U}$.
\begin{prop}\label{prop:transfer maps for davis isotropy}
	The map $d^1_{T,U}$ is non zero only when $U\subset T$ and is given by the following transfer map: 
	\begin{eqnarray*}
		d^1_{T,U}: H_q(W_T; \Z_T)&\to& H_q(W_U; \Z_U).
	\end{eqnarray*}
	On the chain level we compute $H_q(W_T; \Z_T)$ as homology of $\Z_T \otimes_{W_T} F_{W_T}$ for $F_{W_T}$ a projective resolution of $\Z$ over $\Z{W_T}$ and we compute $H_q(W_U; \Z_U)$ as homology of~$\Z_U \otimes_{W_U} F_{W_T}$. Let $m\otimes x$ be in $\Z_T \otimes F_{W_T}$ and $W_U\backslash W_T$ be a set of coset representatives for $W_U$ in $W_T$. Then on the chain level the transfer map is
	\begin{eqnarray*}
		d^1_{T,U}: 
		m \otimes x &\mapsto& \sum_{g \in W_U\backslash W_T} m\cdot g^{-1} \otimes g\cdot x.
	\end{eqnarray*}
	\begin{proof}
		This proof follows the description of the $d^1$ differential for the isotropy spectral sequence in~\cite[VII.8]{Brown}.
		Recall that an orbit representative for a $p$-cell is $eW_T$ for $T$ in $\mathcal{S}$ and $|T|=p$. The set $\mathcal{F}_T$ of cells in the image of the cellular differential~$\cd(W_T)$ is given by cells~$wW_U$ with $|U|=(p-1)$ and $wW_U \subset W_T$. This is satisfied if and only if $U\subset T$ and $w\in W_T$ by Lemma \ref{lem:coset inclusion coxeter}. Since $W_T$ is the stabiliser of the cell $eW_T$, the orbit set~$(\mathcal{F}_T\slash W_T)$ is given by~$\{U \subset T \mid |U|=p-1\}$, which is a subset of~$\mathcal{O}_{p-1}$. The intersection~$\operatorname{Stab}(W_T)\cap \operatorname{Stab}(W_U)=W_T \cap W_U=W_U$ and the action of $W_U$ on $\Z_T$ precisely mimics the action of $W_U$ on $\Z_U$. Therefore it follows that		
		$$
			d^1|_{H_q(W_T; \Z_T)}=\sum_{U \in \mathcal{F}_T/W_T}t_{T, U}
		$$
		
		\noindent where $t_{T,U}$ is the transfer map		
		$
		t_{T,U}:H_q(W_T; \Z_T) \to H_q(W_U;\Z_U).
		$
		 
		Note that cycles in $H_q(W_T; \Z_T)$ are represented by chains in $\Z_T \otimes F_{W_T}$. Letting~$m\otimes x$ be an element on the chain level yields the formula, where the transfer map on the chain level is computed via \cite[III.9]{Brown}.
	\end{proof}
\end{prop}

To compute $H_2(W;\Z)$ and $H_3(W;\Z)$  we consider the~$E^{\infty}$ groups on the~$p+q=2$ (red) diagonal and the~$p+q=3$ (blue) diagonal of Figure \ref{fig:isotropy ss} respectively. Entries on the~$E^1$ page are given by summing over finite Coxeter groups with generating set a certain size, and the classification of finite Coxeter groups from Theorem \ref{prop:classification of finite Coxeter} provides a finite selection of possible groups for each size of generating set.

\begin{lem}\label{lem:inculsion of gps gives inclusion of ss coxeter}
	Given a Coxeter system $(W,S)$, let $V\hookrightarrow W$ be a parabolic subgroup. Then there is a map of isotropy spectral sequences~$
	E(V)\to E(W)
	$ that is an inclusion on the $E^1$ page.
	\begin{proof}
		The inclusion $j:V\hookrightarrow W$ induces an inclusion $W_V\mathcal{S}_V\subset W\mathcal{S}$, hence a map between the realisations $i:\Sigma_V \hookrightarrow \Sigma_W$, and therefore a map~$$C_p(\Sigma_V, \Z)\overset{i_*}{\longrightarrow}C_p(\Sigma_W, \Z).$$ We have the following diagram:		
		
		\centerline{\xymatrix@C=20mm{
				E^1_{p,q}(V)\ar@{=}[d]\ar[r]^{}&\ar@{=}[d]E^1_{p,q}(W)\\
				H_q(V; C_p(\Sigma_V, \Z)) \ar[r]^{H_q(j_*;i_*)}\ar[d]^{\cong}&H_q(W; C_{p}(\Sigma_W, \Z))\ar[d]^{\cong}\\
				\bigoplus_{\substack{U \in \mathcal{S}_V\\ | U | = p}} H_q(W_U; \Z_U)\ar@{.>}[r]^{}& 
				\bigoplus_{\substack{T \in \mathcal{S}\\ | T | = p}} H_q(W_T; \Z_T) \\
		}}
		\noindent where the dotted map is induced by the map on $p$-cells on the central row. Every spherical subgroup of $V$ is also a spherical subgroup of $W$, corresponding to a map between the $p$-cells represented by these spherical subgroups. Therefore the dotted map is an inclusion of summands.
		\noindent Since the $d^1$ differential is defined via the transfer map on each summand, all $d^1$ differentials in $E(V)$ will map under the inclusion to the same differential in $E(W)$. The inclusion on the~$E^1$ page therefore induces a map of spectral sequences on further pages.
	\end{proof}
\end{lem}

\subsection{Pairings on the isotropy spectral sequence}\label{sec:pairings}

We consider a pairing of spectral sequences, for use in Section \ref{sec:h3 further differentials zero}. We follow~\cite{May} and recall Section 4 on products. For filtered complexes $A$, $B$ and $C$, if a pairing~$A \otimes B \to C$ is a morphism of filtered complexes, i.e.~ if $F_pA\cdot F_qB \subset F_{p+q}C$, then this induces a morphism of spectral sequences
$$
E^r(A\otimes B)\to E^r(C).
$$
Combining this with the K\"{u}nneth map $E^r(A)\otimes E^r(B) \to E^r(A\otimes B)$ (which is induced by the K\"{u}nneth map on homology on the $E^1$ page) defines a pairing
$$
\phi:E^r(A)\otimes E^r(B)\to E^r(C)
$$
which satisfies the Leibniz formula for differentials, i.e.~for $x$ in $E^r(A)$ and $y$ in~$E^r(B)$ the pairing satisfies
$$
d_C^r(\phi(x\otimes y))= \phi(d^r_A(x)\otimes y )+ (-1)^{deg(x)}\phi (x \otimes d^r_B(y)).
$$

For finite Coxeter groups~$W_U$ and~$W_V$, let~$W_X=W_U \times W_V$ where~$X:=U \sqcup V$ as in Section \ref{sec:cox products and subgroups}. For the remainder of this section we fix the following notation: let~$W_I$ be the Coxeter group corresponding to~$I\in\{V,U,X\}$. Let~$S_I$ be the generating set of~$W_I$ and let~$\mathcal{S}_I$ be~$\mathcal{S}$ for the Coxeter system $(W_I,I)$ (see Definition~\ref{def:mathcalS}). Let~$\Sigma_I$ be the Davis complex~$\Sigma_{W_I}$ and~$F^I$ be a projective resolution of~$\Z$ over~$\Z{W_I}$. Let~$E(I)$ denote the isotropy spectral sequence for~$W_I$. Then~$E(I)$ is the spectral sequence related to the double complex~$F^I\otimes C(\Sigma_I,\Z)$ (see~\cite[VII.7]{Brown}). Denote the double complex by~$I_{p,q}$ and the associated total complex~$TI$. Then the spectral sequence~$E(I)$ has corresponding filtration $$F_p((TI)_n)=\underset{i\leq p}{\bigoplus} I_{n-i,i}.$$

\begin{lem}\label{lem:pairing on davis}
	The product map $W_U\times W_V \to W_X$ determines a map on chain complexes 
	$$C_i(\Sigma_U, \Z)\otimes C_j(\Sigma_V,\Z) \to C_{i+j}(\Sigma_{X},\Z).$$
	\begin{proof}
		The product map induces a map of posets
		\begin{eqnarray*}
			W_U\mathcal{S}_U \times W_V\mathcal{S}_V &\to& W_X \mathcal{S}_X\\
			(uW_{T_U}, vW_{T_V}) &\mapsto& uv(W_{T_U \sqcup T_V}).
		\end{eqnarray*}
		\noindent
		This in turn induces a map on their realisations~$
		\Sigma_U\times \Sigma_V \to \Sigma_X,
		$	 
		\noindent which is the map giving the decomposition $\Sigma_X=\Sigma_U\times \Sigma_V$ in Lemma \ref{lem:product of Davis}.
		\noindent Consider $C_i(\Sigma_I, \Z)$ and note that $p$-cells of $\Sigma_I$ are represented by cosets $wW_T$ where $T \in \mathcal{S}_I$. Given an~$i$-cell of $\Sigma_U$ represented by $uW_{T_1}$ and a $j$-cell of $\Sigma_V$ represented by $vW_{T_2}$ we use the above poset map and define an $(i+j)$-cell of $\Sigma_X$ represented by $uvW_{T_1\sqcup T_2}$. This gives a pairing $C_i(\Sigma_U, \Z)\otimes C_j(\Sigma_V,\Z) \to C_{i+j}(\Sigma_{X},\Z)$.		
	\end{proof}
\end{lem}

\begin{prop}\label{prop:pairing}
	The map 
	$$
	\Phi:E^r(U)\otimes E^r(V) \to E^r(X)
	$$
	induced by the pairings $F^U_k \otimes F^V_l \to F^X_{k+l}$ and $C_i(\Sigma_U, \Z)\otimes C_j(\Sigma_V,\Z) \to C_{i+j}(\Sigma_{X},\Z)$, gives a pairing of spectral sequences, under which the differentials satisfy the Leibniz formula. 
	
		\begin{proof}
		We apply the hypothesis of~\cite[Section~4]{May} and show that the map
		$
		TU \otimes TV \to TX
		$
		is a morphism of filtered complexes. We have on the $n$th level that 
		$$F_p((TI)_n)=\bigoplus_{i\leq p} I_{n-i,i}=\bigoplus_{i\leq p}F^I_{n-i}\otimes C_i(\Sigma_I, \Z)
		$$
		for $I$ in $\{U,V,X\}$. Since $W_U \times W_V=W_X$, there is a pairing $F^U_k \otimes F^V_l \to F^X_{k+l}$ (e.g.~$F^X=F^U\otimes F^V$ \cite[V.1.1]{Brown}). Putting this together with the pairing $C_i(\Sigma_U, \Z)\otimes C_j(\Sigma_V,\Z) \to C_{i+j}(\Sigma_{X},\Z)$ from Lemma~\ref{lem:pairing on davis} gives 
		$$
		F_p(TU)\cdot F_q(TV)\subset F_{p+q}(TX)
		$$
		as required in \cite{May}.
		
	\end{proof}
\end{prop}
	
	\begin{thm}\label{lem:pairing May}
		
		Under the decomposition on the $E^1$ page of the spectral sequence
		$$
		E^1_{p,q}(I)=H_q(F^I_*\otimes_{W_I} C_p(\Sigma_I,\Z))\cong\bigoplus_{\substack{\bar{I} \in \mathcal{S}_I\\ | \bar{I} | = p}} H_q(W_{\bar{I}}; \Z_{\bar{I}}) 
		$$	
		\noindent the pairing of Proposition~\ref{prop:pairing} induces a pairing $\Phi_*$, given by the K\"{u}nneth map when restricted to individual summands
		$$
		\Phi_*:H_q(W_{\bar{U}}; \Z_{\bar{U}})\otimes  H_{q'}(W_{\bar{V}}; \Z_{\bar{V}}) \overset{\times}{\to}H_{q+q'}(W_{\bar{U}}\times W_{\bar{V}};  \Z_{\bar{U}}\otimes \Z_{\bar{V}})\overset{\cong}{\to}H_{q+q'}(W_{\bar{X}}; \Z_{\bar{X}}).
		$$
		\noindent It follows that the differentials in the isotropy spectral sequence for the Davis complex satisfy a Leibniz formula with respect to the pairing $\Phi_*$.
		
		\begin{proof}
				We now consider this pairing under the decomposition on the $E^1$ page of the isotropy spectral sequence:
			$$
			E^1_{p,q}(I)=H_q(F^I_*\otimes_{W_I} C_p(\Sigma_I,\Z))\cong\bigoplus_{\substack{\bar{I} \in \mathcal{S}_I\\ | \bar{I} | = p}} H_q(W_{\bar{I}}; \Z_{\bar{I}}) 
			$$	
			\noindent described in~\cite{Brown}[VII]. Under this decomposition the above isomorphism restricted to a single summand on the right is given by the following map $\iota_*$, induced by inclusion $\iota:	F_*^T\otimes_{W_T} C_p(\Sigma_T, \Z_T)\to F^W_*\otimes_W C_p(\Sigma_W,\Z)$.
			\ms
			
			\centerline{\xymatrix@R=5mm{
					H_q(F_*^T\otimes_{W_T} C_p(\Sigma_T, \Z_T))\ar[r]^{\iota_*}&H_q(F^W_*\otimes_W C_p(\Sigma_W,\Z))\\
					H_q(W_T; \Z_T) \ar[r]\ar@{=}[u]& \ar@{=}[u]H_q(F^W_*\otimes_W C_p(\Sigma_W,\Z)).				
			}}
			\ms
			\noindent If a Coxeter group $W_X$ arises as a product $W_X=W_U\times W_V$, then the pairing $\Phi$, along with the~$E^1$ decomposition for each group gives the following diagram.
			
			\centerline{\xymatrix@R=10mm{
					\big(E_{p,q}^1(U) \ar@{=}[d]\ar@{}@<-1.5ex>[r]^{\bigotimes}& E_{p',q'}^1(V)\big)\ar@{=}[d]\ar[r]^-{\Phi}&\ar@{=}[d]E^1_{(p+p'),(q+q')}(X)\\	
					\big(H_q(F^U_*\otimes_{W_U} C_p(\Sigma_U,\Z))\ar@{}@<-1.5ex>[r]^{\bigotimes}& H_{q'}(F^V_*\otimes_{W_V} C_{p'}(\Sigma_V,\Z))\big)\ar[r]^-{\Phi}	& H_{q+q'}(F^X_*\otimes_{W_X} C_{p+p'}(\Sigma_X,\Z))\\
					\big(\underset{\substack{\bar{U} \in \mathcal{S}_U\\ | \bar{U} | = p}}{\bigoplus} H_q(W_{\bar{U}}; \Z_{\bar{U}})\ar@{}@<-1.5ex>[r]^{\bigotimes}\ar[u]_{\cong}^{\oplus \iota_*} &\ar[u]_{\cong}^{\oplus \iota_*}\underset{\substack{\bar{V} \in \mathcal{S}_V\\ | \bar{V} | = p'}}{\bigoplus} H_{q'}(W_{\bar{V}}; \Z_{\bar{V}})\big) \ar@{.>}[r]^-{\Phi_*}&\underset{\substack{\bar{X} \in \mathcal{S}_X\\ | \bar{X} | = p+p'}}{\bigoplus} H_{q+q'}(W_{\bar{X}}; \Z_{\bar{X}}). \ar[u]_{\cong}^{\oplus \iota_*}\\
			}}
			\ms
			
			\noindent The isomorphisms are induced by the component-wise inclusions given by $\iota_*$ on each summand. The map $\Phi_*$ is defined such that the diagram commutes, i.e.~it is induced by $\Phi$ and the two vertical isomorphisms. On each summand of the bottom left factor it is given by the composite
			$$
			H_q(W_{\bar{U}}; \Z_{\bar{U}})\otimes  H_{q'}(W_{\bar{V}}; \Z_{\bar{V}}) \overset{\times}{\to}H_{q+q'}(W_{\bar{U}}\times W_{\bar{V}};  \Z_{\bar{U}}\otimes \Z_{\bar{V}})\overset{\cong}{\to}H_{q+q'}(W_{\bar{X}}; \Z_{\bar{X}})
			$$
			where here $\bar{X}:=\bar{U}\sqcup\bar{V}$. Here the first map is given by the homology cross product \cite[V.3]{Brown}, and the second map is due the fact that if  $W_{\bar{U}}\times W_{\bar{V}}=W_{\bar{X}}$ then the orientation modules satisfy $\Z_{\bar{U}}\otimes \Z_{\bar{V}}\cong \Z_{\bar{X}}$. This map is precisely the K\"{u}nneth map on homology.
			Extending this component-wise definition to the tensor product of the summations gives the map $\Phi_*$ that lifts to the map $\Phi$ on the top row. 
		\end{proof}
	\end{thm}

\section{Calculation of $H_2(W;\Z)$}\label{sec:h2 calc coxeter}

From Section \ref{sec:isotropy spectral seq for davis}, the isotropy spectral sequence for $(W,S)$ has $E^1$ page as in Figure \ref{fig:isotropy ss}, and the $E^{\infty}$ page will give filtration quotients of $H_2(W;\Z)$ on the $p+q=2$ diagonal. We compute the diagonal on the $E^2$ page and note that no further differentials affect this diagonal, so the result follows.

It the following let~$(W,S)$ be a Coxeter system and~$E^1_{p,q}:=E^1_{p,q}(W)$ be the~$E^1$ terms of the isotropy spectral sequence for the Davis complex of~$(W,S)$. 

\begin{prop}\label{prop:zeros for h2 diagonal} The terms~$E^1_{0,2}$ and~$E^1_{1,1}$ are zero.
\begin{proof}
	We have~$E^1_{0,2}=H_2(W_\emptyset; \Z_\emptyset)=0$, since~$W_{\emptyset}$ is the trivial group. The $E^1_{1,1}$ term is given by 
	$$E^1_{1,1}=\underset{t \in S}{\oplus} H_1(W_t; \Z_t),$$ where the non-trivial group element $t$ acts by negation. Then $H_1(W_t; \Z_t)=0$ follows from taking the standard projective resolution for a cyclic group of order 2 and these coefficients.
\end{proof}
\end{prop}

\subsection{ Homology at $E^1_{2,0}$}\label{sec:sequence}
Recall that
$$E^1_{2,0}=\underset{| T | =2}{\underset{T \in \mathcal{S}}{\oplus}}   H_0(W_T; \Z_T) .$$
\noindent From Proposition~\ref{prop:zeros for h2 diagonal}, this will be the only contributing group to the~$p+q=2$ diagonal on the $E^\infty$ page. We start by computing $E^2_{2,0}$, which is given by the homology of the following sequence

\centerline{\xymatrix{
		\underset{t \in \mathcal{S}}{\oplus} H_0(W_t; \Z_t) &\ar[l]^{d^1}
		\underset{| T | =2}{\underset{T \in \mathcal{S}}{\oplus}}   H_0(W_T; \Z_T) & \ar[l]^{d^1}\underset{| T | =3}{\underset{T \in \mathcal{S}}{\oplus}}   H_0(W_T; \Z_T). 
}}

\noindent Recall that the $d^1$ differential is given by the transfer map defined in Proposition~\ref{prop:transfer maps for davis isotropy}.

\begin{lem} For all $T$ in $\mathcal{S}$ such that $\vert T \vert >0$,~$H_0(W_T; \Z_T)=\Z_2.$
	\begin{proof}
		The zeroth homology is given by the coinvariants of the coefficient module $\Z_T$ under the group action. Since in our case each group generator acts as multiplication by $-1$ we compute homology to be the group $\Z_2$.
	\end{proof}
\end{lem}

For $X \in \mathcal{S}$, let $1_X$ be the generator for the summand $H_0(W_X, \Z_X)$ of $E^1_{p,0}$.

\begin{lem}\label{lem:transfer on bottom row}
	When $U$ is a subset of $T$, the transfer map for the bottom row of the spectral sequence is
	\begin{eqnarray*}
		d^1_{T,U}:H_0(W_T; \Z_T)=\Z_2 &\to& H_0(W_{U}; \Z_{U})=\Z_2\\
		1_T &\mapsto& 
		\left\{
		\begin{array}{ll}
			0  & \text{ if  }| W_T | / | W_{U}|\text{ is even }  \\
			1_{U}  &\text{ if  }| W_T | / | W_{U}|\text{ is odd.}
		\end{array}
		\right.
	\end{eqnarray*} 
	\begin{proof}
	  
	From~\cite[III.9.(B)]{Brown}, the transfer map acts on coinvariants as
		\begin{eqnarray*}
			d^1_{T,U}:H_0(W_T; \Z_T)=\Z_2 &\to& H_0(W_{U}; \Z_{U})=\Z_2\\
			1_T &\mapsto& \sum_{g \in W_{U}\backslash W_T} g\cdot 1_{U}=\sum_{g \in W_{U}\backslash W_T} 1_{U}
		\end{eqnarray*}
		\noindent since $g\cdot 1=\pm1$ is in the class of $1$ in $\Z_{U}\slash W_{U}$. Noting that we are mapping into~$\Z_2$ and the number of entries in the sum is $| W_T | / | W_{U}|$ completes the proof.
	\end{proof}
\end{lem}

\begin{lem}\label{lem:h2 transfer 2 to 1}
	When $U$ has cardinality 1 and $T=\{s,t\}$ has cardinality 2 the transfer map $d^1$ restricted to the $T$ summand is given by
	$$
	d^1|_{H_0(W_{T};\Z_{T})}(1_T)=
	\left\{
	\begin{array}{ll}
	1_s+1_t  & \mbox{if } m(s,t) \text{ odd} \\
	0 & \mbox{if } m(s,t) \text{ even}.
	\end{array}
	\right.
	$$
	\begin{proof}
		Note that $| W_x |= 2$ for all $x$ in $S$ and since $W_{\{s,t\}}$ is isomorphic to a dihedral group, $| W_{\{s,t\}}|=2 \times m(s,t)$. Apply Lemma~\ref{lem:transfer on bottom row} to compute the differential.
	\end{proof}
\end{lem}

\begin{defn}
	We say that a Coxeter group with generating set $T=\{s,t,u\}$  is \emph{of type X} if if $W_T=W({\bf I}_2(p))\times W({\bf A}_1)$ and $p\geq 3$ is odd, i.e.~$\D_{W_T}$ has the form:
	
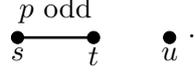
\begin{center}
		\begin{tikzpicture}[scale=0.2, baseline=0]
		\draw[fill= black] (5,0) circle (0.4) node[below] {$s$};
		\draw (7.5,1.8) node {$p$ odd};
		\draw[line width=1] (5,0) -- (10,0);			
		\draw[fill= black] (10,0) circle (0.4) node[below] {$t$};	
		\draw[fill= black] (15,0) circle (0.4) node[below] {$u$};
		\end{tikzpicture}.
	\end{center}
\end{defn}

\begin{lem}\label{lem:h2 transfer 3 to 2}\label{coset quoients}
	If $T=\{s,t,u\}$ then $d^1$ restricted to the~$H_0(W_{T};\Z_{T})$ summand is
	$$
	d^1|_{H_0(W_{T};\Z_{T})}(1_T)=
	\left\{
	\begin{array}{ll}
	1_{\{s,u\}}+1_{\{t,u\}}  & \mbox{if }  W_T \text{ is of type } X\\
	0 & \mbox{otherwise.}
	\end{array}
	\right.
	$$
	\begin{proof}
		There are a finite number of Coxeter diagrams that may represent $W_T$, given by Theorem \ref{prop:classification of finite Coxeter}. The order of these groups and their rank two subgroups is documented in the table below, where $p \geq 2$.
		\begin{center}
			\begin{tabular}[center]{|c|c||c|c|c|c|} 
				\hline
				$W_T$&$\D_W$ & $|W_T|$ & $|W_{\{s,t\}}|$&  $|W_{\{s,u\}}|$  &  $|W_{\{t,u\}}|$ \\
				\hline &&&&&\\[-1em]
				$W({\bf A}_3)$&\begin{tikzpicture}[scale=0.15, baseline=0]
				\draw[fill= black] (5,0) circle (0.5) node[below] {$s$};
				\draw[line width=1] (5,0) -- (15,0);
				\draw[fill= black] (10,0) circle (0.5) node[below] {$t$};
				\draw[fill= black] (15,0) circle (0.5)node[below] {$u$};
				\end{tikzpicture}  & 24 &6 & 4& 6\\	\hline
				$W({\bf D}_3)$ &\begin{tikzpicture}[scale=0.15, baseline=0]
				\draw[fill= black] (5,0) circle (0.5) node[below] {$s$};
				\draw (7.5,1.5) node {$4$};
				\draw[line width=1] (5,0) -- (15,0);
				\draw[fill= black] (10,0) circle (0.5) node[below] 	{$t$};
				\draw[fill= black] (15,0) circle (0.5)node[below] {$u$};
				\end{tikzpicture}  &48 &8 &4 &6\\	\hline
				$W({\bf H}_3)$&
				\begin{tikzpicture}[scale=0.15, baseline=0]
				\draw[fill= black] (5,0) circle (0.5) node[below] {$s$};
				\draw (7.5,1.5) node {$5$};
				\draw[line width=1] (5,0) -- (15,0);
				\draw[fill= black] (10,0) circle (0.5) node[below] 	{$t$};
				\draw[fill= black] (15,0) circle (0.5)node[below] {$u$};
				\end{tikzpicture}  & 120& 10& 4&6\\	\hline
				$W({\bf I}_2(p))\times W({\bf A}_1)$&\begin{tikzpicture}[scale=0.15, baseline=0]
				\draw[fill= black] (5,0) circle (0.5) node[below] {$s$};
				\draw (7.5,1.5) node {$p$};
				\draw[line width=1] (5,0) -- (10,0);
				\draw[fill= black] (10,0) circle (0.5) node[below] 	{$t$};
				\draw[fill= black] (15,0) circle (0.5)node[below] {$u$}; 
				\end{tikzpicture}  & 4$p$& 2$p$ & 4 &  4\\
				\hline
			\end{tabular}
		\end{center}
		Calculating $| W_T | / | W_{T'}|$ for $T'\subset T$ in each of these cases and applying Lemma~\ref{lem:transfer on bottom row} completes the proof.
	\end{proof}
\end{lem}

\begin{prop}\label{prop:h2 splitting calc}
	The homology at $E^1_{2,0}$ is given by 	
	$$H_0(\D_{\bullet \bullet}; \Z_2) \oplus \Z_2[E(\D_{\operatorname{even}})] \oplus H_1(\D_{\operatorname{odd}};\Z_2)$$
	\noindent where the diagrams are as defined in Definition \ref{def:diagrams for h2} and are viewed as 1-dimensional complexes.
	\begin{proof} Consider the calculations of the transfer maps in Lemmas \ref{lem:h2 transfer 2 to 1} and \ref{lem:h2 transfer 3 to 2}, and observe the following splitting:

		\centerline{\xymatrix{
				\underset{t \in {S}}{\oplus} H_0(W_t; \Z_t)  \ar@{=}[ddd] &\ar@{=}[d]\underset{| T | =2}{\underset{T \in \mathcal{S}}{\oplus}}   H_0(W_T; \Z_T)\ar[l]^{d^1} &\underset{| T | =3}{\underset{T \in \mathcal{S}}{\oplus}}   H_0(W_T; \Z_T)\ar[l]^{d^1} \\
				&\ar@{}@<-1.5ex>[d]^{\bigoplus} \underset{ m(s,t)=2}{\underset{T =\{s,t\}}{\oplus}}   H_0(W_T; \Z_T)&	\ar[l]^{d^1}\underset{W_T\text{ type }X
				}{\oplus}H_0(W_T;\Z_T)\ar@{^{(}->}[u]\\
				&\ar@{}@<-1.5ex>[d]^{\bigoplus}\underset{ m(s,t)\neq 2\text{ even}}{\underset{T =\{s,t\}}{\oplus}}   H_0(W_T; \Z_T)&\\
				\underset{t \in {S}}{\oplus} H_0(W_t; \Z_t) &\underset{ m(s,t)\text{ odd}}{\underset{T =\{s,t\}}{\oplus}}   H_0(W_t; \Z_t).\ar[l]^-{d^1}&
		}}
		\ms
		
		\noindent Calculating the homology of the top row in turn gives a splitting
		\begin{eqnarray*}
			&\operatorname{coker}\Big(\underset{W_T\text{ type }X
			}{\oplus}H_0(W_T,\Z_T)\overset{d^1}{\to} \underset{ m(s,t)=2}{\underset{T =\{s,t\}}{\oplus}}   H_0(W_T; \Z_T)\Big)&\\
			&\bigoplus&\\
			&\underset{ m(s,t)\neq 2\text{ even}}{\underset{T =\{s,t\}}{\oplus}}   H_0(W_T; \Z_T)&\\
			&\bigoplus&\\
			&\ker\Big(\underset{ m(s,t)\text{ odd}}{\underset{T =\{s,t\}}{\oplus}}   H_0(W_T; \Z_T)\overset{d^1}{\to} \underset{t \in \mathcal{S}}{\oplus} H_0(W_T; \Z_T)\Big).&\\
		\end{eqnarray*}
		
		\noindent We now define an isomorphism $\varepsilon=\varepsilon_1 \oplus \varepsilon_2 \oplus \varepsilon_3$ from these three groups, to the three groups in the statement of the proposition:
		$$
		H_0(\D_{\bullet \bullet}; \Z_2) \oplus \Z_2[E(\D_{\operatorname{even}})] \oplus H_1(\D_{\operatorname{odd}};\Z_2).
		$$
		
		The map between the first groups is as follows:
		\begin{eqnarray*}
			\varepsilon_1:\operatorname{coker}\Big(\underset{W_T\text{ type }X
			}{\oplus}H_0(W_T,\Z_T)\overset{d^1}{\to} \underset{ m(s,t)=2}{\underset{T =\{s,t\}}{\oplus}}   H_0(W_T; \Z_T)\Big) &\to& H_0(\D_{\bt \bt}; \Z_2)\\
			1_{\{s,t\}} &\mapsto& [\{s,t\}],
		\end{eqnarray*}
		where $[\{s,t\}]$ is the generator for the summand of $H_0(\D_{\bullet \bullet}; \Z_2)$ corresponding to the connected component containing $\{s,t\}$.
		
		Recall from Lemma \ref{lem:h2 transfer 3 to 2} that the transfer map on summands $H_0(W_{\{s,t,u\}};\Z_T)$ is given by $d^1(1_{\{s,t,u\}})= 1_{\{s,u\}}+1_{\{t,u\}}$ if $W_T$ is of type X. Therefore generators of~$H_0(W_T;\Z_T)$ for triples of type $X$ get mapped to sums of generators of~$H_0(W_T;\Z_T)$ corresponding to commuting pairs. These are exactly vertices of~$\D_{\bt \bt}$, and a triple of type $X$ gives the corresponding edge of~$\D_{\bt \bt}$. Therefore the map $\varepsilon_1$ is well defined and moreover it is an isomorphism.
		
		For $\Z_2[E(\D_{\operatorname{even}})]$, let $\{s,t\}$ be the basis element corresponding to the edge between $s$ and $t$, and note that edges only exist if $m(s,t)$ is even and greater than 2. Then $\varepsilon_2$ is the isomorphism defined by
		\begin{eqnarray*}
			\varepsilon_2: \underset{ m(s,t)\neq 2,\text{even}}{\underset{T =\{s,t\}}{\oplus}}   H_0(W_T; \Z_T)&\to&\Z_2[E(\D_{\operatorname{even}})]\\
			1_{\{s,t\}} &\mapsto& \{s,t\}.
		\end{eqnarray*}
		
		For $H_1(\D_{\operatorname{odd}};\Z_2)$, note that $\D_{\operatorname{odd}}$ has no 2-cells, so  $ H_1(\D_{\operatorname{odd}};\Z_2)$ is the kernel of the cellular differential~$\cd:C_1\to C_0$, where~$C_1=~\Z_2[E(\D_{\operatorname{odd}})]$,~$C_0=\Z_2[S]$, and $\cd(\{s,t\})=s+t$. Recall from Lemma \ref{lem:h2 transfer 2 to 1} that the transfer map is given on summands $H_0(W_{\{s,t\}};\Z_T)=\Z_2$ by $d^1(1_{\{s,t\}})= 1_s+1_t$ if $m(s,t)$ is odd. Therefore we define a chain map:
		\begin{eqnarray*}
			\underset{ m(s,t)\text{ odd}}{\underset{T =\{s,t\}}{\oplus}}   H_0(W_T; \Z_T)&\to&\Z_2[E(\D_{\operatorname{odd}})]\\
			1_{\{s,t\}} &\mapsto& \{s,t\}
		\end{eqnarray*}		
		\noindent and this map induces an isomorphism~$	\varepsilon_3$ between homologies.
	\end{proof}
\end{prop}

\subsection{Proof of Theorem \ref{thm:INTRO THM A}}\label{sec:h2 proof of theorem}

\begin{thm} Given a finite rank Coxeter group $W$ with diagram $\D_W$, recall from Definition \ref{def:diagrams for h2} the definition of the diagrams $\D_{\bt \bt}$, $\D_{\operatorname{odd}}$ and $\D_{\operatorname{even}}$. Then there is a natural isomorphism
	$$H_2(W;\Z)=	H_0(\D_{\bullet \bullet}; \Z_2) \oplus \Z_2[E(\D_{\operatorname{even}})] \oplus H_1(\D_{\operatorname{odd}};\Z_2)$$
	\noindent where in the first and final term of the right-hand side the diagrams are viewed as cell complexes.
	\begin{proof}
		The $p+q=2$ diagonal of the isotropy spectral sequence in Figure \ref{fig:isotropy ss} gives filtration quotients of $H_2(W;\Z)$ on the $E^\infty$ page. The $E^2$ page has only one non-zero term on this diagonal:
		$$
		E^2_{2,0}H_0(\D_{\bullet \bullet}; \Z_2) \oplus \Z_2[E(\D_{\operatorname{even}})] \oplus H_1(\D_{\operatorname{odd}};\Z_2).
		$$
		All differentials $d^r$ for $r\geq 2$ with source or target the $E_{2,0}$ position either originate at, or map to a zero group. Therefore the $p+q=2$ diagonal on the limiting $E^\infty$ page is given by the diagonal on the $E^2$ page. Since there is only one non zero group on the diagonal, there are no extension problems and this group gives $H_2(W;\Z)$ as required.
	\end{proof}
\end{thm}

\section{Calculation of $H_3(W;\Z)$}\label{SEC:H3 FOR COXETER}

Recall the isotropy spectral sequence for the Coxeter group $W$ has $E^1$ page as shown in Figure \ref{fig:isotropy ss} in Section \ref{sec:isotropy spectral seq for davis}, and the $E^{\infty}$ page gives $H_3(W;\Z)$ (up to extension) on the $p+q=3$ diagonal.

In Section \ref{sec:resolution for finite Coxeter} the free resolution for finite Coxeter groups by De Concini and Salvetti \cite{DCSal} is introduced and the chain map between resolutions is computed in Section \ref{sec:collapse}. Using these tools, we compute the $E^2$ page of the spectral sequence on the~$p+q=3$ diagonal. Following this, Section \ref{sec:h3 further differentials zero} proves that all further differentials to and from this diagonal are zero. The possible extension problems arising on the limiting page at this diagonal are treated and discussed in Section~\ref{sec:extension problem} and all of these computations are fed into the proof of Theorem \ref{thm:INTRO THM B} in Section \ref{sec:h3 proof of theorem}.

\subsection{Free resolution for Coxeter groups}\label{sec:resolution for finite Coxeter}

In~\cite{DCSal}, De Concini and Salvetti introduce a free resolution of~$\Z$ over~$\Z W$ for a finite Coxeter group $W$. We use this throughout this section to calculate the low dimensional homologies of finite Coxeter groups that appear as summands in the $E^1$ entries of the spectral sequence.

\begin{defn}\label{def:DCSal chain complex}
	Let~$(W,S)$ be a Coxeter system for a finite Coxeter group~$W$. Let~$(C_*, \delta_*)$ be the chain complex with~$C_k$ the free~$\Z W$ module with basis elements~$e(\G)$. Here~$\G$ is a flag of subsets of the generating set~$S$ with cardinality~$k$, that is~$\G$ in~$S_k$ where:
$$
S_k:= \{\G=(\G_1\supset \G_2\supset \cdots ) \mid \G_1 \subset S, \sum_{i\geq 1} |\G_i| =k\}.
$$
For $\tau$ in $\G_i$, let $W_{\G_i}^{\Gamma_i\backslash \{\tau\} }$ be the set of minimal left coset representatives of $W_{\G_i\backslash \{\tau\} }$ in $W_{\G_i}$. Then $\delta_k:C_k\to C_{k-1}$ is $\Z W$-linear and defined on basis elements by
\begin{equation}\label{eq:free res differential}
\delta_k e(\G)= \sum_{\substack{i\geq 1 \\ |\G_i|> |\G_{i+1}|}} \sum_{\tau \in \G_{i}} \sum_{\substack{\beta \in W_{\G_i}^{\Gamma_i\backslash \{\tau\} } \\ \beta^{-1} \G_{i+1}\beta \subset \G_i\backslash \{\tau \}}} (-1)^{\alpha(\G, i, \tau, \beta)}\beta e(\G')
\end{equation}

\noindent where the flag $\G'$ in $C_{k-1}$ is given by
$$
\G':=(\G_1\supset \cdots \supset \G_{i-1}\supset (\G_i\backslash \{\tau\}) \supset \beta^{-1} \G_{i+1} \beta \supset\beta^{-1} \G_{i+2} \beta \supset \cdots)
$$

\noindent and the exponent $\alpha(\G, i, \tau, \beta)$ is defined below. The differential is well defined from Lemma \ref{lem:shortest coset element coxeter}. We choose an ordering for the set of generators $S$ and let $\sigma(\beta, \G_k)$ be the number of inversions, with respect to this ordering, in the map $\G_k \to \beta^{-1}\G_k \beta$. We let $\mu(\G_i,\tau)$ be the number of generators in $\G_i$ which are less than or equal to $\tau$ in the ordering on $S$. Then the exponent is described by the following formula:
\begin{equation*}
\alpha(\G, i, \tau, \beta)=i \cdot \lng(\beta) + \sum_{k=1}^{i-1}|\G_k| + \mu(\G_i, \tau)+\sum_{k=i+1}^d \sigma(\beta, \G_k).
\end{equation*}
\end{defn}

\noindent During this proof we adopt the convention that the generators are always ordered alphabetically (e.g.~$s<t<u$). We also denote the generator corresponding to a flag of length $d$, $(\G_1\supset \G_2\supset \cdots \supset\G_d)$ by $\G_{\G_1\supset \G_2\supset \cdots \supset \G_d}$, where we omit the set notation for each $\G_i$. For example $\G_s$, $\G_{s\supset s}$, or $\G_{ s,t \supset s}$ (which corresponds to $\G=\{s,t\}\supset \{s\}$).

\begin{thm}[{\cite{DCSal}}]
	The chain complex $(C_*,\delta_*)$ from Definition \ref{def:DCSal chain complex} is a free resolution of $W$ over $\Z W$.
\end{thm}

%

\begin{exam}\label{ex:c2 resolution}
	We give an example of the resolution for finite Coxeter groups with one generator $S=\{s\}$, from $C_3$ to $C_0$.
	
	\ms
	
	\centerline{\xymatrix@R=1mm@C=17mm {
			C_3=\langle\G_{s \supset s \supset s}\rangle \ar[r]^{\delta_3=(s-1)}				& 	C_2=\langle\G_{s \supset s}\rangle \ar[r]^{\delta_2=(1+s)}				& C_1=\langle \G_s\rangle  \ar[r]^{\delta_1=(s-1)} 	& C_0=\langle\G_\emptyset \rangle 	
	}}

	\ms
	
	\noindent The differential from $\G_s$ to $\G_\emptyset$ is given by the following formula, noting that coset representatives of $W_\emptyset$ in $W_s$ are $e$ and $s$. We recall the formula for $\delta_k(e(\G))$ from Equation \eqref{eq:free res differential}.
	\begin{eqnarray*}
		\delta_1 (\G_s)&=&\sum_{\beta=e, s } (-1)^{\alpha(\G_s, 1, s, \beta)}\beta \G_\emptyset 
		= (s-1)\G_\emptyset
	\end{eqnarray*}	
	\noindent where we compute	
	\begin{eqnarray*}
		\alpha(\G_s, 1, s, e)&=& 1 \lng(e) + \sum_{k=1}^{0}|\G_k| + \mu(s, s)= 0 + 0 + 1 = 1 \\
		\alpha(\G_s, 1, s, s)&=& 1 \lng(s) + \sum_{k=1}^{0}|\G_k| + \mu(s, s)=
		 1 + 0 + 1 = 2.
	\end{eqnarray*}	
	\noindent Similarly the differential $\delta_2:C_2 \to C_1$ is given by	
	\begin{eqnarray*}
		\delta_2 (\G_{s \supset s})&=&  \sum_{\beta=e, s } (-1)^{\alpha(\G_{s \supset s}, 2, s, \beta)}\beta \G_s 
		= (1+s)\G_s
	\end{eqnarray*}	
	\noindent where we compute 
	\begin{eqnarray*}
		\alpha(\G_{s \supset s}, 2, s, e)&=& 2 \lng(e) + \sum_{k=1}^{1}|\G_k| + \mu(s, s)=0 + 1 + 1 =2 \\
		\alpha(\G_{s \supset s}, 2, s, s)&=& 2 \lng(s) + \sum_{k=1}^{1}|\G_k| + \mu(s, s)=2 + 1 + 1 = 4.
	\end{eqnarray*}	
	\noindent Finally, the differential $\delta_3:C_3 \to C_2$ is given by 
	\begin{eqnarray*}
		\delta_3 (\G_{s\supset s \supset s})
		=\sum_{\beta=e, s } (-1)^{\alpha(\G_{s\supset s \supset s}, 3, s, \beta)}\beta \G_{s\supset s } 
		= (s-1)\G_{s\supset s }
	\end{eqnarray*}	
	\noindent where we compute 
	\begin{eqnarray*}
		\alpha(\G_{s\supset s \supset s}, 3, s, e)&=& 3 \lng(e) + \sum_{k=1}^{2}|\G_k| + \mu(s, s)= 0 + 2 + 1  = 3 
\\
		\alpha(\G_{s\supset s \supset s}, 3, s, s)&=& 3 \lng(s) + \sum_{k=1}^{2}|\G_k| + \mu(s, s)= 3 + 2 + 1  = 6.
	\end{eqnarray*}
\end{exam}

\begin{defn}\label{def:p(s,t;m)}
	Define $p(s,t;j)$ to be the alternating product of $s$ and $t$ of length $j$, \emph{ending} in an $s$ (as opposed to $\altprod(s,t;j)$ which is the alternating product \emph{starting} in an $s$) i.e. \vspace{-.35cm}
	$$
	p(s,t;j)=\overbrace{\ldots sts}^\text{length j}.
	$$
\end{defn}

\begin{exam}\label{ex:res for dihedral}
	Consider the resolution for finite Coxeter groups with two generators $S=\{s, t\}$, from $C_3$ to $C_0$ and with $m(s,t)$ finite. Then the formulas for differentials which do not follow from the previous example are:
	\begin{eqnarray*}
				\delta_2(\G_{s, t}) &=& {\displaystyle \sum_{j=0}^{m(s,t)-1} (-1)^{j+1} p(s,t;j)\G_t + \sum_{g=0}^{m(s,t)-1} (-1)^{g+2}p(t,s;g) \G_s}	 \\
				\delta_3(\G_{s, t \supset s})  &=&{\begin{cases}
						(1-p(t,s;m(s,t)-1))\G_{s\supset s}-(1+s)\G_{st}&
						\text{if } m(s,t) \text{ even}\\ \G_{s\supset s}-p(s,t;m(s,t)-1)\G_{t\supset t}-(1+s)\G_{st}&\text{if } m(s,t) \text{ odd}
				\end{cases}}
				\\
			\delta_3(\G_{s, t \supset t}) &=&{\begin{cases} (-1+p(s,t;m(s,t)-1))\G_{t\supset t}-(1+t)\G_{st}&\text{if } m(s,t) \text{ even }\\
						-\G_{t\supset t}+p(t,s;m(s,t)-1)\G_{s\supset s}-(1+t)\G_{st}&\text{if } m(s,t) \text{ odd }
				\end{cases}}
	\end{eqnarray*}
\end{exam}

Recall we wish to compute homologies of finite Coxeter groups $W_T$ with twisted coefficients $\Z_T$, in which the action of the generators on $\Z_T$ is given by negation. To calculate the twisted homologies we tensor the resolution with $\Z$ under the group action. We show this in the case of our two examples.

\begin{exam}\label{ex: C_2 twisted res}
	We consider the resolution of Example \ref{ex:c2 resolution} tensored with $\Z$ under the group action:
	\ms
	
	\centerline{\xymatrix@R=3mm {
			{\begin{array}{rl} &\Z \underset{W_s}{\otimes} C_3\\=&\langle 1 \otimes\G_{s \supset s \supset s\rangle}\end{array}} \ar[r]^-{\delta_3=(-2)} &	{\begin{array}{rl}&\Z \underset{W_s}{\otimes} C_2\\=&\langle1 \otimes\G_{s \supset s}\rangle \end{array}} \ar[r]^-{\delta_2=(0)} 			& {\begin{array}{rl}&\Z \underset{W_s}{\otimes}C_1\\=&\langle 1 \otimes \G_s \rangle \end{array}} \ar[r]^-{\delta_1=(-2)} 	&{\begin{array}{rl}& \Z \underset{W_s}{\otimes}C_0\\=&\langle 1 \otimes \G_\emptyset \rangle. \end{array}}
	}}
\ms
\noindent Here the differentials are calculated as follows:
\begin{eqnarray*}
\delta_3(1 \otimes \G_{s \supset s\supset s})&=&1\otimes( (s-1)\G_{s\supset s })=-2(1\otimes \G_{s\supset s})\\
\delta_2(1 \otimes \G_{s \supset s})&=& 1\otimes((1+s)\G_s)=0\\
\delta_1(1 \otimes \G_{s})&=& 1\otimes ((s-1)\G_\emptyset)=-2(1\otimes \G_\emptyset).
\end{eqnarray*}
\end{exam}

\begin{exam} \label{ex: dihedral twisted res}
	We consider the computations of differentials in Example \ref{ex:res for dihedral} and upon tensoring with $\Z$ under the group action, this gives the following differentials:	
			\begin{eqnarray*}
					\delta_2(1 \otimes \G_{s, t}) &=&	 -m(s,t)(1 \otimes \G_t)+m(s,t)(1\otimes \G_s)\\				
				\delta_3(	1 \otimes \G_{s, t \supset s}) &= &  {\begin{cases}
							2(1 \otimes \G_{s\supset s})&\text{ if } m(s,t) \text{ even }\\
							1 \otimes\G_{s\supset s}-1 \otimes \G_{t\supset t}&\text{ if } m(s,t) \text{ odd }
					\end{cases}}
					\\
					\delta_3(1 \otimes \G_{s, t \supset t}) &= &  {\begin{cases}-2(1 \otimes \G_{t\supset t})&\text{ if } m(s,t) \text{ even }\\
							-1 \otimes\G_{t\supset t}+1 \otimes\G_{s\supset s}&\text{ if } m(s,t) \text{ odd. }		
					\end{cases}}			
			\end{eqnarray*}
	
\end{exam}

\subsection{Collapse map}\label{sec:collapse}
In this section we define a chain map, which we call the \emph{collapse map}, between De Concini and Salvetti's resolution for a finite Coxeter group $W$, and for a subgroup $W_T$ \cite{DCSal}.

Recall that in the isotropy spectral sequence for the Davis complex, the $d^1$ differential has the form of a transfer map, given in Proposition \ref{prop:transfer maps for davis isotropy}. In the following sections we calculate these twisted homology groups using the De Concini and Salvetti resolution. Upon applying the transfer map to a generator of~$H_*(W_T;\Z_T)$, the image will be in terms of the resolution for the group $W_T$. However we require the image to be in terms of the resolution for $W_U$ and so we apply the collapse map in the appropriate degree to achieve this.

We first recall the following lemmas from \cite{GeckPfeiffer}. Recall from Definition \ref{def: pi(s,t;m)} that~$\altprod(a,b ; k)$ is defined to be the word of length $k$, given by the alternating product of $a$ and $b$.

\begin{lem}[{Deodhar's Lemma, \cite[2.1.2]{GeckPfeiffer} }]\label{lem:deodhar}
	For $(W,S)$ a Coxeter system, let $W_T$ be a spherical subgroup of a finite Coxeter group $W$, let $v$ be $(T, \emptyset)$-reduced (Definition \ref{def:(-,-) reduced coxeter}) and let $s$ be in $S$. Then either $vs$ is $(T, \emptyset)$-reduced or $vs=tv$ for some $t$ in $T$.
\end{lem}

\begin{lem}[{\cite[1.2.1]{GeckPfeiffer}}] \label{lem:delta(s,t) at end}
	If $s,u$ are in $S$, $m(s,u)$ is finite, and $w$ in $W$ satisfies~$\lng(ws)<\lng(w)$ and $\lng(wu)<\lng(w)$ then it follows $w=w'(\altprod(s,u;m(s,u)))$ where~$w'$ is $(\emptyset, W_{\{s,u\}})$-reduced.
\end{lem}

\begin{defn}\label{def:collapse map}
	Denote the De Concini - Salvetti resolution for $(W,S)$ by $(C_*, \delta_*)$ and for the subgroup $(W_T,T)$ by $(D_*, \delta_*)$. We define the \emph{collapse map in degree} $i$ to be the $W_T$-equivariant linear map $f_i:C_i \to D_i$ for $0\leq i \leq 2$ as shown below.	
	\centerline{\xymatrix@R=10mm {
			\ar[r]^{\delta_3} &C_2 \ar[r]^{\delta_2} 	\ar[d]_{f_2}			& C_1 \ar[r]^{\delta_1} \ar[d]_{f_1}	& C_0	\ar[r]^{\delta_0}\ar[d]_{f_0}	& \Z \ar@{=}[d]	\\
			\ar[r]^{\delta_3} &D_2 \ar[r]^{\delta_2} 				& D_1 \ar[r]^{\delta_1} 	& D_0	\ar[r]^{\delta_0}	& \Z.	
	}}
	
	\noindent As a $\Z[W]$ module, $C_*$ has basis given by $e(\G)$, so as a $\Z[W_T]$ module, $C_*$ has basis given by $v\cdot e(\G)$, for $v$ a $(T,\emptyset)$-reduced element of $W$. We therefore define~$f_i$ on~$v\cdot e(\G)$ and extend the map linearly and $W_T$-equivariantly. By Lemma \ref{lem:deodhar} for~$s~\in~S$,~$vs$ is either~$(T, \emptyset)$-reduced or~$vs=tv$ for some~$t$ in~$T$. This gives us the cases in each definition.	
	$$f_0 (v\G_\emptyset)=\G_\emptyset,$$ 	
	$$
	f_1(v\G_s)=
	\begin{cases} 
	0 & vs \text{ is } (T, \emptyset) \text{ reduced} \\
	\G_t & vs= tv \text{ for } t \in T \\
	\end{cases}
	$$
	$$
	f_2(v\G_{s \supset s})=
	\begin{cases} 
	0 & vs \text{ is } (T, \emptyset) \text{ reduced} \\
	\G_{t \supset t} & vs= tv \text{ for } t \in T \\
	\end{cases}
	$$
	$$
	f_2(v\G_{su})= \begin{cases} 
	\G_{tr} & vs=tv \text{ and } vu=rv \text{ for } t,r \in T \\
	0 & \text{ otherwise. }
	\end{cases}
	$$ 
	
	
\end{defn}
The remainder of this section is devoted to proving that $f_*$ is a chain map.

\begin{lem}\label{lem:f0 colllapse commutes}
	The following square commutes:
	
	\centerline{\xymatrix@R=10mm {
			C_0 \ar[r]^{\delta_0} \ar[d]_{f_0}	& \Z\ar@{=}[d]		\\
			D_0 \ar[r]^{\delta_0} 	& \Z.
	}}
	
	\begin{proof}
		Let $w$ in $W$. For each basis element $w\G_\emptyset$, the square is given by
		
		\centerline{
			\xymatrix@R=10mm {
				w\G_\emptyset \ar[r]^{\delta_0} \ar@{|->}[d]_{f_0}	& \Z\ar@{=}[d]		\\
				f_0 (w\G_\emptyset) \ar[r]^{\delta_0} 	& \Z.		
		}}
		
		\noindent Since $f_0$ is defined $W_T$-equivariantly, if $w=w'v$ for $w'$ in $W_T$ and $v$ a $(T, \emptyset)$-reduced element then from Definition \ref{def:collapse map} 
		$$
		f_0 (w\G_\emptyset)=f_0 (tv\G_\emptyset)=t\cdot f_0 (v\G_\emptyset)=t\G_\emptyset.
		$$  
		
		\noindent It follows since $\delta_0$ maps all generators to $1$ that the square commutes.
	\end{proof}
\end{lem}

\begin{lem}\label{lem:f1 colllapse commutes}
	The following square commutes
	
	\centerline{\xymatrix@R=10mm {
			C_1 \ar[r]^{\delta_1} \ar[d]_{f_1}	& C_0\ar[d]_{f_0}		\\
			D_1 \ar[r]^{\delta_1} 	& D_0.	
	}}
	
	\begin{proof}
		Since all maps are $W_T$-equivariant, we need only consider the square on generators multiplied by a~$(T,\emptyset)$-reduced element~$v$. We recall the image of $\delta_1$ from Example \ref{ex:c2 resolution}.
		
		\centerline{\xymatrix@R=10mm {
				v\G_s \ar[r]^-{\delta_1} \ar@{|->}[d]_{f_1}	&v(s-1)\G_\emptyset\ar@{|->}[d]_{f_0}		\\
				f_1(v\G_s) \ar[r]^-{\delta_1} 	& f_0(v(s-1)\G_\emptyset).\\
		}}
		
		\noindent Here thetwo cases for the element $vs$, given by Lemma \ref{lem:deodhar}, give the following cases for $f_0$, from Definition \ref{def:collapse map}:		
		$$
		f_0(v(s-1)\G_\emptyset)=
		\begin{cases} 
		0 & vs \, (T, \emptyset) \text{ reduced} \\
		(t-1)\G_\emptyset & vs= tv. \\
		\end{cases}
		$$
		
		\noindent This is precisely the image of $f_1(v\G_s)$ from Definition \ref{def:collapse map}, under the differential~$\delta_1$. Therefore the square commutes.
	\end{proof}
\end{lem}

	For $s$ and $u$ in $S$, consider the following three cases, given by Lemma~\ref{lem:deodhar}:
	\begin{enumerate}
		\item  Neither $vs$ or $vu$ are $(T, \emptyset)$-reduced, that is $vs=tv$ and $vu=rv$ for $t$ and $r$ in $T$.
		\item One of $vs$ and $vu$ is $(T, \emptyset)$-reduced, without loss of generality let $vs=tv$ and $vu$ is $(T, \emptyset)$-reduced.
		\item Both $vs$ and $vu$ are $(T, \emptyset)$-reduced.
	\end{enumerate}
	Recall from Definition \ref{def:p(s,t;m)} that $p(s,u;m)$ is the alternating product of $s$ and $u$ of length $m$ ending in~$s$.
	\begin{lem} \label{lem:degree 2 su} 
    We have that
	$$	\begin{array}{c}
	f_1\left(v\big(\sum_{j=0}^{m(s,u)-1} (-1)^{j+1} p(s,u;j)\G_u +  \sum_{g=0}^{m(s,u)-1} (-1)^{g+2}p(u,s;g) \G_s\big)\right)\\=
	\begin{cases} 
	\delta_2(\G_{tr}) &\text{ in Case } (1) \\
	0 &\text{ in Case } (2) \\
	0 &\text{ in Case }(3).
	\end{cases}
	\end{array}	$$
	\begin{proof}
		For Case (1), since $f_1$ acts $W_T$-equivariantly,
		$$
		f_1(v (p(s,u;j)\G_u))=f_1(p(t,r;j)v\G_u)=p(t,r;j)(f_1(v\G_u))=p(t,r;j)\G_r
		$$
		and similarly 
		$
		f_1(vp(u,s;g) \G_s)=p(r,t;g)\G_t.
		$
		Furthermore, $m(t,r)=m(s,u)$ since 
		\begin{eqnarray*}
			\altprod(t,r;m(s,u))v=v\altprod(s,u;m(s,u))=v\altprod(u,s;m(s,u))=\altprod(r,t;m(s,u))v,
		\end{eqnarray*}
		
		\noindent and by right multiplication by $v^{-1}$, $\altprod(t,r;m(s,u))=\altprod(r,t;m(s,u))$, so~$m(t,r)$ is a divisor of  $m(s,u)$. Applying a similar argument in reverse gives~$m(s,u)$ is a divisor of~$m(t,r)$, and so $m(s,u)=m(t,r)$. Therefore since $f_1$ acts linearly, we have
		\begin{eqnarray*}
			&&f_1\left(v\big(\sum_{j=0}^{m(s,u)-1} (-1)^{j+1} p(s,u;j)\G_u +  \sum_{g=0}^{m(s,u)-1} (-1)^{g+2}p(u,s;g) \G_s\big)\right)\\
			&=& \sum_{j=0}^{m(t,r)-1} (-1)^{j+1} p(t,r;j)\G_r +  \sum_{g=0}^{m(t,r)-1} (-1)^{g+2}p(r,t;g)\G_t 
			=\delta_2(\G_{tr}).
		\end{eqnarray*}
		
		For Case (2), we first prove that if~$vs=tv$ and~$vu$ is~$(T, \emptyset)$-reduced, it follows~$v(\altprod(u,s;k))$ is also~$(T, \emptyset)$-reduced for all~$2\leq k \leq m(s,u)-1$. Note that since~$vs=tv$, from Lemma~\ref{lem:deodhar}~$\lng(vs)>\lng(v)$. Suppose~$v(\altprod(u,s;k))$ is not~$(T, \emptyset)$-reduced and choose minimal~$k$ for which this is the case. Then for some~$q$ in~$T$ it follows~$v(\altprod(u,s;k))=qv(\altprod(u,s;k-1))$ and so~$w=v(\altprod(u,s;k))$ satisfies the hypothesis of Lemma~\ref{lem:delta(s,t) at end}, that is~$\lng(wu)<\lng(w)$ and~$\lng(ws)<\lng(w)$. Therefore $$w=w'\altprod(u,s;m(s,u)))=v(\altprod(u,s;k)).$$
		By right multiplication by~$(\altprod(u,s;k))^{-1}$ we have~$v=w' p(s,u;m(s,u)-k)$. Therefore $v$ satisfies~$\lng(vs)<\lng(v)$, but this contradicts $vs=tv$. Therefore $v(\altprod(u,s;k))$ is also $(T, \emptyset)$-reduced for all~$2\leq~k~\leq~m(s,u)-1$. Computing $f_1$ it follows:
		$$
		f_1(v (p(s,u;j)\G_u))= \begin{cases} 
		f_1(v(\altprod(u,s;j)\G_u))=0 & j \text{ is even, }j\neq m(s,u)-1 \\
		t\cdot f_1(v\altprod(u,s;j-1)\G_u)=t\cdot 0=0 & j \text{ is odd, }j\neq m(s,u)-1  \\
		f_1(v \altprod(u,s;m(s,t)-1)\G_u)=\G_t & j= m(s,u)-1  \text{ and is even}\\
		t\cdot f_1(v \altprod(u,s;m(s,t)-2)\G_u)=t\cdot0 & j= m(s,u)-1  \text{ and is odd}
		\end{cases}
		$$
		and similarly 
		$$
		f_1(vp(u,s;g) \G_s)=\begin{cases} 
		f_1(v\G_s)=\G_t & g=0 \\
		t\cdot f_1(v\altprod(u,s;g-1) \G_s)=t\cdot 0=0  & g\text{ is even, } g\notin\{0,m(s,u)-1\}\\
		f_1(v\altprod(u,s;g) \G_s)=0& g \text{ is odd, }g\neq m(s,u)-1\\		
		t\cdot f_1(v \altprod(u,s;m(s,t)-2) \G_s)=t\cdot 0=0&g= m(s,u)-1  \text{ and is even}\\
		f_1(v \altprod(u,s;m(s,t)-1) \G_s)= \G_t& g= m(s,u)-1  \text{ and is odd}
		\end{cases}
		$$
		so it follows in the setting of Case (2) we have
		\begin{eqnarray*}
			&&f_1\left(v\big(\sum_{j=0}^{m(s,u)-1} (-1)^{j+1} p(s,u;j)\G_u +  \sum_{g=0}^{m(s,u)-1} (-1)^{g+2}p(u,s;g) \G_s\big)\right)\\
			&=& \begin{cases} 
				\G_t + (-1)^{m(s,t)-1+2}\G_t =0 & \text{ if } m(s,u) \text{ even }\\ 
				
				\G_t +(-1)^{m(s,u)-1+1}\G_t =0 & \text{ if } m(s,u) \text{ odd.}
			\end{cases}
		\end{eqnarray*}
		
		For Case (3), if both $vs$ and $vu$ are $(T, \emptyset)$-reduced, by the same argument as in Case (2), $v(\altprod(u,s;k))$ and $v(\altprod(s,u;k))$ are also $(T, \emptyset)$-reduced for~$2\leq k \leq m(s,u)$. Computing $f_1$ 
		in the setting of Case (3) gives:
		$$
		f_1\left(v\big(\sum_{j=0}^{m(s,u)-1} (-1)^{j+1} p(s,u;j)\G_u +  \sum_{g=0}^{m(s,u)-1} (-1)^{g+2}p(u,s;g) \G_s\big)\right)=0.
		$$
	\end{proof}
\end{lem}

\begin{lem}\label{lem:f2 colllapse commutes}
	The following square commutes
	
	\centerline{\xymatrix@R=10mm {
			C_2 \ar[r]^{\delta_2} \ar[d]_{f_2}	& C_1\ar[d]_{f_1}		\\
			D_2 \ar[r]^{\delta_2} 	& D_1.
	}}
	\begin{proof}
		
		Since all maps are $W_T$-equivariant, let $v$ be a $(T, \emptyset)$-reduced element and consider the square on generators left-multiplied by $v$. We recall the image of~$\delta_2$ from Example \ref{ex:res for dihedral}. We must consider both forms of generators of $C_2$:
		
		\centerline{	
			\xymatrix@R=10mm {
				v\G_{s\supset s} \ar[r]^(0.4){\delta_2} \ar@{|->}[d]_{f_2}	&v(1+s)\G_s\ar@{|->}[d]_{f_1}		\\
				f_2(v\G_{s\supset s}) \ar[r]^(0.4){\delta_2} 	& f_1(v(1+s)\G_s)\\
			}\hspace{1cm}
			\xymatrix@R=10mm {
				v\G_{s,u} \ar[r]^(0.25){\delta_2} \ar@{|->}[d]_{f_2}	&{\begin{array}{l}v\big(\sum_{j=0}^{m(s,u)-1} (-1)^{j+1} p(s,t;j)\G_u \\+  \sum_{g=0}^{m(s,u)-1} (-1)^{g+2}p(u,s;g) \G_s\big)\end{array}}	\ar@{|->}[d]_{f_1}	\\
				f_2(v\G_{s,u}) \ar[r]^(0.25){\delta_2} 	&{\begin{array}{l} f_1\big(v\big(\sum_{j=0}^{m(s,u)-1} (-1)^{j+1} p(s,u;j)\G_t \\+  \sum_{g=0}^{m(s,u)-1} (-1)^{g+2}p(u,s;g) \G_s\big)\big).\end{array}}
		}}
		
		\noindent Computing $f_1(v(1+s)\G_s)$ we have
		$$
		f_1(v(1+s)\G_s)=
		\begin{cases} 
		0 & vs \text{ is } (T, \emptyset) \text{ reduced} \\
		(1+t)\G_t & vs= tv.\\
		\end{cases}
		$$
		
		\noindent This is precisely the image of $f_2(v\G_{s \supset s})$ from Definition \ref{def:collapse map}, under the differential~$\delta_2$. Therefore the left hand square commutes.
		
		The bottom right entry of the right hand square is given in Lemma \ref{lem:degree 2 su}. This is precisely the image of $f_2(v\G_{s,u})$ from Definition \ref{def:collapse map}, under the differential $\delta_2$. Therefore the right hand square commutes.
	\end{proof}
\end{lem}

\begin{prop}
	The maps $f_0$, $f_1$ and $f_2$ in Definition \ref{def:collapse map} form part of a chain map $f_\bt:C_\bt \to D_\bt$.
	\begin{proof}
		This is a consequence of Lemmas \ref{lem:f0 colllapse commutes}, \ref{lem:f1 colllapse commutes} and \ref{lem:f2 colllapse commutes},
		which show that all the required squares commute.
	\end{proof}
\end{prop}

In the following sections the tools we have developed are utilised to compute the~$E^2$ terms of the isotropy spectral sequence for the Davis complex. When a proof is omitted, this is due to the fact it is a straightforward calculation of homology. All omitted proofs can be found in~\cite[Appendix B]{Boydthesis}.

\begin{lem}
	For~$r\geq1$, we have $E^1_{0,r}=H_r(W_\emptyset; \Z_\emptyset)=0$.
\end{lem}

It follows that the~$E^1_{0,3}$ term of the diagonal is zero on the $E^{\infty}$ page.
	
\subsection{Homology at $E^1_{1,2}$}

We use the De Concini - Salvetti resolution \cite{DCSal} and the transfer (Proposition \ref{prop:transfer maps for davis isotropy}) and collapse (Definition \ref{def:collapse map}) maps to compute the differentials for the following section of the spectral sequence:

\ms
\centerline{\xymatrix@R=3mm {
		0=H_2(W_\emptyset; \Z_\emptyset)	& \ar[l]_{d^1} \color{teal}{\underset{t \in S}{\oplus} H_2(W_t; \Z_t)}  &  \ar[l]_{d^1} \underset{| T | =2}{\underset{T \in \mathcal{S}}{\oplus}}  H_2(W_T; \Z_T). }}

Let~$W_t$ and~$W_T$ be as in the above sequence, and~$T=\{s,t\}$
\begin{lem}\label{lem:h2 of 1 generator group}\label{lem:h2 2 generators}
	In terms of the De Concini - Salvetti resolution, the homologies in the above sequence are
	$H_2(W_t; \Z_t)= \Z_2$, generated by $1 \otimes \G_{t \supset t} $, and
	$$
	H_2(W_T; \Z_T)=\begin{cases}
	\Z_2 \oplus \Z_2 & \text{ if } m(s,t) \text{ is even}\\	 
	\Z_2  & \text{ if } m(s,t) \text{ is odd},
	\end{cases}	
	$$
	generated by $1 \otimes\G_{s\supset s}$ and $1 \otimes \G_{t\supset t}$ when $m(s,t)$ is even, with these generators being identified when $m(s,t)$ is odd.
\end{lem}

\begin{lem} \label{lem:degree 2 map st to s}
	For~$u$ in~$T$,~$d^1_{T,u}$ is given by
	\begin{eqnarray*}
		d^1_{T,u}:H_2(W_{\{s,t\}}; \Z_T) &\to& H_2(W_u; \Z_u) \\
		1\otimes \G_{s\supset s} &\mapsto&  1\otimes \G_{u\supset u}
	\end{eqnarray*}
	if $m(s,t)$ is odd, and the zero map if~$m(s,t)$ is even.
	\begin{proof}
		We apply the transfer map from Proposition \ref{prop:transfer maps for davis isotropy} to the generator(s) of $H_2(W_{\{s,t\}}; \Z_T)$ followed by the degree two collapse map $f_2$ from Definition \ref{def:collapse map}.
	\end{proof}
\end{lem}

\begin{prop}
	The $E^2_{1,2}$ entry of the isotropy spectral sequence for $(W,S)$ is given by $H_0(\D_{\operatorname{odd}};\Z_2)$.
	\begin{proof}
		On the $E^1$ page we compute homology of the sequence
		
		\centerline{\xymatrix@R=3mm {
				0	& \ar[l]_{d^1} \underset{t \in S}{\oplus} \Z_2  &  \ar[l]_(0.7){d^1}  \underset{m(s,t) \text{ even}}{\underset{T=\{s,t\}}{\underset{T \in \mathcal{S}}{\oplus}}} (\Z_2 \oplus \Z_2 ) \underset{m(s,t) \text{ odd}}{\underset{T=\{s,t\}}{\underset{T \in \mathcal{S}}{\oplus}}} \Z_2
	}}
		
		\noindent The left hand map is the zero map and the right hand map is defined via Lemma~\ref{lem:degree 2 map st to s}. Applying the splitting technique as in the proof of the $H_2(W;\Z)$ calculation (Proposition~\ref{prop:h2 splitting calc}), gives homology equal to $H_0(\D_{\operatorname{odd}};\Z_2)$ as required.
	\end{proof}
\end{prop}

\subsection{Homology at $E^1_{2,1}$}
The $E^1$ page at $E^1_{2,1}$ has the following form:
\ms

\centerline{\xymatrix@R=3mm {
		\underset{t \in S}{\oplus} H_1(W_t; \Z_t) & \ar[l]_{d^1} \color{teal}{\underset{| T | =2}{\underset{T \in \mathcal{S}}{\oplus}}  H_1(W_T; \Z_T)}&  \ar[l]_{d^1} \underset{| T | =3}{\underset{T \in \mathcal{S}}{\oplus}} H_1(W_T; \Z_T). }}

\begin{prop}\label{prop:h1 of 3 gen groups}
	The first homology $H_1(W_T; \Z_T)$ is as follows for finite $W_T$ with $T=\{s,t,u\}$. Generators are given by the De Concini - Salvetti resolution for $W_T$: we set $$\alpha=(1 \otimes \G_s)-(1\otimes \G_t) \text{ and }\beta=(1 \otimes \G_s)-(1\otimes \G_u).$$	
	\begin{center}
		\begin{tabular}[center]{|c|c|c|c|}
			\hline
			$W_T$ & $\D_{W_T}$ & $ H_1(W_T; \Z_T) $& Generator \\
			\hline &&&\\[-1em]
			$W({\bf A}_3)$&\begin{tikzpicture}[scale=0.15, baseline=0]
			\draw[fill= black] (5,0) circle (0.5) node[below] {$s$};
			\draw[line width=1] (5,0) -- (15,0);
			\draw[fill= black] (10,0) circle (0.5) node[below] {$t$};
			\draw[fill= black] (15,0) circle (0.5)node[below] {$u$}; 
			\end{tikzpicture}  & $\Z_3$ & $\alpha$\\	\hline
			$W({\bf B}_3)$& \begin{tikzpicture}[scale=0.15, baseline=0]
			\draw[fill= black] (5,0) circle (0.5) node[below] {$s$};
			\draw (7.5,1.2) node {4};
			\draw[line width=1] (5,0) -- (15,0);
			\draw[fill= black] (10,0) circle (0.5) node[below] 	{$t$};
			\draw[fill= black] (15,0) circle (0.5)node[below] {$u$}; 
			\end{tikzpicture}  & $\Z_2$ & $\alpha=\beta$\\	\hline
			$W({\bf H}_3)$&\begin{tikzpicture}[scale=0.15, baseline=0]
			\draw[fill= black] (5,0) circle (0.5) node[below] {$s$};
			\draw (7.5,1.2) node {5};
			\draw[line width=1] (5,0) -- (15,0);
			\draw[fill= black] (10,0) circle (0.5) node[below] 	{$t$};
			\draw[fill= black] (15,0) circle (0.5)node[below] {$u$};
			\end{tikzpicture}  & $0$&\\	\hline
			${\begin{array}{c c} \\ W({\bf I}_2(p)) \times W({\bf A}_1)\\p\geq2\end{array}}$ & \begin{tikzpicture}[scale=0.15, baseline=0]
			\draw[fill= black] (5,0) circle (0.5) node[below] {$s$};
			\draw (7.5,1.2) node {p};
			\draw[line width=1] (5,0) -- (10,0);
			\draw[fill= black] (10,0) circle (0.5) node[below] 	{$t$};
			\draw[fill= black] (15,0) circle (0.5)node[below] {$u$};
			\end{tikzpicture}  &${\begin{array}{c r}
				\\\Z_2 \oplus \Z_2  &\text{ if } p \text{ is even}\\	 
				\Z_2   &\text{ if } p \text{ is odd}\\\end{array}}$&${\begin{array}{c r}
				\\ \alpha, \beta  &\text{ if } p \text{ is even}\\	 
				\beta  &\text{ if } p \text{ is odd}\end{array}}$ \\
			& & & \\
			\hline
		\end{tabular}
	\end{center}
\end{prop}

\begin{prop}\label{prop:h1 of 2 gen group}
	When $T=\{s,t\}$, $H_1(W_T;\Z_T)=\Z_{m(s,t)}$ with generator in the De Concini - Salvetti resolution given by $\gamma=1 \otimes \G_s - 1 \otimes \G_t.$
\end{prop}

\begin{prop}
	Let~$s$ in~$S$. Then~$H_1(W_s;\Z_s)=0$.
\end{prop}

We now introduce some notation. If $H_i(W_T; \Z_T)$ only has one generator, then we represent that generator in the~$E^1_{p,q}$ summation of homologies by drawing the diagram $\D_{W_T}$. We represent~$d^1|_{H_i(W_T; \Z_T)}$ by drawing a map from the diagram~$\D_{W_T}$ to the diagrams representing generators in the image of~$d^1|_{H_i(W_T; \Z_T)}$, with signs and scalar multiplication as required. In some cases $H_i(W_T; \Z_T)$ has either zero or two generators, but in these cases there are non non-zero differentials.

%

\begin{prop}\label{prop:d1 E^2 for H_1s}
	The non-zero differentials on the $E^1$ page at $E^1_{2,1}$ are given as in the diagram below.
	
	\centerline{\xymatrix@R=0mm {
			\underset{t \in S}{\oplus} H_1(W_t; \Z_t) &  \underset{| T | =2}{\underset{T \in \mathcal{S}}{\oplus}}  H_1(W_T; \Z_T)\ar[l]_{d^1} &  \underset{| T | =3}{\underset{T \in \mathcal{S}}{\oplus}} H_1(W_T; \Z_T)\ar[l]_{d^1} \\
			&	\begin{tikzpicture}[scale=0.15, baseline=0]
			\draw[fill= black] (5,0) circle (0.5) node[below] {$s$};
			\draw[line width=1] (5,0) -- (10,0);
			\draw[fill= black] (10,0) circle (0.5) node[below] {$t$};
			\end{tikzpicture} \ominus	\begin{tikzpicture}[scale=0.15, baseline=0]
			\draw[line width=1] (10,0) -- (15,0);
			\draw[fill= black] (10,0) circle (0.5) node[below] {$t$};
			\draw[fill= black] (15,0) circle (0.5)node[below] {$u$}; 
			\end{tikzpicture}&\begin{tikzpicture}[scale=0.15, baseline=0]
			\draw[fill= black] (5,0) circle (0.5) node[below] {$s$};
			\draw[line width=1] (5,0) -- (15,0);
			\draw[fill= black] (10,0) circle (0.5) node[below] {$t$};
			\draw[fill= black] (15,0) circle (0.5)node[below] {$u$}; 
			\end{tikzpicture}\ar@{|->}[l]\\
			&
			\begin{tikzpicture}[scale=0.15, baseline=0]
			\draw[fill= black] (8,0) circle (0.5) node[below] {$s$};
			\draw[fill= black] (15,0) circle (0.5)node[below] {$u$};
			\end{tikzpicture} \oplus
			\begin{tikzpicture}[scale=0.15, baseline=0]
			\draw[fill= black] (10,0) circle (0.5) node[below] 	{$t$};
			\draw[fill= black] (15,0) circle (0.5)node[below] {$u$};
			\end{tikzpicture}&	\begin{tikzpicture}[scale=0.15, baseline=0]
			\draw[fill= black] (5,0) circle (0.5) node[below] {$s$};
			\draw (7.5,1.4) node {p \text{ odd }};
			\draw[line width=1] (5,0) -- (10,0);
			\draw[fill= black] (10,0) circle (0.5) node[below] 	{$t$};
			\draw[fill= black] (15,0) circle (0.5)node[below] {$u$};
			\end{tikzpicture} \ar@{|->}[l]
	}}
	
	\begin{proof}
		This proof involves calculating the differential $d^1$ via the transfer and collapse maps. This can be calculated by hand, but we use Python and the PyCox package \cite{Geck}. These calculations can be found in \cite[Appendix B]{Boydthesis}.
	\end{proof}		
	
\end{prop}

\begin{prop}
	Recall from Definition \ref{def:diagrams for h3} the diagrams $\D_{\bt \bt}$ and $\D_{{\bf A}_2}$. Then the $E^2_{2,1}$ entry of the isotropy spectral sequence for $(W,S)$ is given by
	$$
	H_0(\D_{\bullet \bullet };\Z_2) \oplus H_0(\D_{{\bf A}_2};\Z_3) \oplus\Big( \underset{m(s,t)>3,\neq \infty}{\oplus}\Z_{m(s,t)}\Big).
	$$
	\begin{proof}
		Consider the $d^1$ differentials at $E^2_{2,1}$, given in Proposition \ref{prop:d1 E^2 for H_1s}, and apply the splitting technique as in Proposition \ref{prop:h2 splitting calc}.
	\end{proof}
\end{prop}

\subsection{Homology at $E^1_{3,0}$}
\begin{lem}\label{lem:d1 at (3,0) for h3}
The non-zero $d^1$ differentials at $E^1_{3,0}$ are given by the following maps
	
	\centerline{	\xymatrix@R=-1mm {
			\underset{| T | =2}{\underset{T \in \mathcal{S}}{\oplus}}  H_0(W_T; \Z_T) & \ar[l]_{d^1} \color{teal}{\underset{| T | =3}{\underset{T \in \mathcal{S}}{\oplus}}  H_0(W_T; \Z_T)}&  \ar[l]_{d^1} \underset{| T | =4}{\underset{T \in \mathcal{S}}{\oplus}}  H_0(W_T; \Z_T) \\
			\underset{| T | =2}{\underset{T \in \mathcal{S}}{\oplus}}  \Z_2 & \ar[l]_{d^1} \underset{| T | =3}{\underset{T \in \mathcal{S}}{\oplus}}  \Z_2&  \ar[l]_{d^1} \underset{| T | =4}{\underset{T \in \mathcal{S}}{\oplus}}   \Z_2\\
			\begin{tikzpicture}[scale=0.15, baseline=0] \draw[fill= black] (10,0) circle (0.5)node[below] 	{$t$};
			\draw[fill= black] (15,0) circle (0.5)node[below] {$u$};
			\end{tikzpicture} + \begin{tikzpicture}[scale=0.15, baseline=0]
			\draw[fill= black] (5,0) circle (0.5) node[below] {$s$};
			\draw[fill= black] (12,0) circle (0.5)node[below] {$u$};
			\end{tikzpicture}
			&  	\ar@{|->}[l]\begin{tikzpicture}[scale=0.15, baseline=0]
			\draw[fill= black] (5,0) circle (0.5) node[below] {$s$};
			\draw (7.5,1.4) node {p odd};
			\draw[line width=1] (5,0) -- (10,0);
			\draw[fill= black] (10,0) circle (0.5) node[below] 	{$t$};
			\draw[fill= black] (15,0) circle (0.5)node[below] {$u$};
			\end{tikzpicture}\\	
			& 	\begin{tikzpicture}[scale=0.15, baseline=0]
			\draw[fill= black] (5,0) circle (0.5) node[below] {$s$};
			\draw[line width=1] (5,0) -- (15,0);
			\draw[fill= black] (10,0) circle (0.5) node[below] {$t$};
			\draw[fill= black] (15,0) circle (0.5)node[below] {$u$}; 
			\end{tikzpicture} + 	\begin{tikzpicture}[scale=0.15, baseline=0]
			\draw[line width=1] (10,0) -- (20,0);
			\draw[fill= black] (10,0) circle (0.5) node[below] {$t$};
			\draw[fill= black] (15,0) circle (0.5)node[below] {$u$};
			\draw[fill= black] (20,0) circle (0.5)node[below] {$v$}; 
			\end{tikzpicture}& \ar@{|->}[l]	\begin{tikzpicture}[scale=0.15, baseline=0]
			\draw[fill= black] (5,0) circle (0.5) node[below] {$s$};
			\draw[line width=1] (5,0) -- (20,0);
			\draw[fill= black] (10,0) circle (0.5) node[below] {$t$};
			\draw[fill= black] (15,0) circle (0.5)node[below] {$u$};
			\draw[fill= black] (20,0) circle (0.5)node[below] {$v$}; 
			\end{tikzpicture}\\
			&	\begin{tikzpicture}[scale=0.15, baseline=0]
			\draw (17.5,1.4) node {q even};
			\draw[line width=1] (15,0) -- (20,0);
			\draw[fill= black] (10,0) circle (0.5) node[below] 	{$t$};
			\draw[fill= black] (15,0) circle (0.5)node[below] {$u$};
			\draw[fill= black] (20,0) circle (0.5)node[below] {$v$}; 	\end{tikzpicture}+ 	
			\begin{tikzpicture}[scale=0.15, baseline=0]
			\draw[fill= black] (7.5,0) circle (0.5) node[below] {$s$};
			\draw (17.5,1.4) node {q even};
			\draw[line width=1] (15,0) -- (20,0);
			\draw[fill= black] (15,0) circle (0.5)node[below] {$u$};
			\draw[fill= black] (20,0) circle (0.5)node[below] {$v$}; 
			\end{tikzpicture} & \ar@{|->}[l]	\begin{tikzpicture}[scale=0.15, baseline=0]
			\draw[fill= black] (5,0) circle (0.5) node[below] {$s$};
			\draw (7.5,1.4) node {p odd};
			\draw (17.5,1.4) node {q even};
			\draw[line width=1] (5,0) -- (10,0);
			\draw[line width=1] (15,0) -- (20,0);
			\draw[fill= black] (10,0) circle (0.5) node[below] 	{$t$};
			\draw[fill= black] (15,0) circle (0.5)node[below] {$u$};
			\draw[fill= black] (20,0) circle (0.5)node[below] {$v$}; 
			\end{tikzpicture}\\
			& { \left( \begin{array}{cccc}	\begin{tikzpicture}[scale=0.15, baseline=0]
					\draw (17.5,1.4) node {q odd};
					\draw[line width=1] (15,0) -- (20,0);
					\draw[fill= black] (10,0) circle (0.5) node[below] 	{$t$};
					\draw[fill= black] (15,0) circle (0.5)node[below] {$u$};
					\draw[fill= black] (20,0) circle (0.5)node[below] {$v$}; 	\end{tikzpicture}&+ &
					\begin{tikzpicture}[scale=0.15, baseline=0]
					\draw[fill= black] (7.5,0) circle (0.5) node[below] {$s$};
					\draw (17.5,1.4) node {q odd};
					\draw[line width=1] (15,0) -- (20,0);
					\draw[fill= black] (15,0) circle (0.5)node[below] {$u$};
					\draw[fill= black] (20,0) circle (0.5)node[below] {$v$}; 
					\end{tikzpicture} &+\\\begin{tikzpicture}[scale=0.15, baseline=0]
					\draw[fill= black] (5,0) circle (0.5) node[below]{$s$};
					\draw (7.5,1.4) node {p odd};
					\draw[line width=1] (5,0) -- (10,0);
					\draw[fill= black] (10,0) circle (0.5) node[below] 	{$t$};
					\draw[fill= black] (17.5,0) circle (0.5)node[below] {$v$}; 
					\end{tikzpicture}&+& 	\begin{tikzpicture}[scale=0.15, baseline=0]
					\draw[fill= black] (5,0) circle (0.5) node[below]{$s$};
					\draw (7.5,1.4) node {p odd};
					\draw[line width=1] (5,0) -- (10,0);
					\draw[fill= black] (10,0) circle (0.5) node[below] 	{$t$};
					\draw[fill= black] (15,0) circle (0.5)node[below] {$u$};
					\end{tikzpicture}\end{array} \right) }	& \ar@{|->}[l]	\begin{tikzpicture}[scale=0.15, baseline=0]
			\draw[fill= black] (5,0) circle (0.5) node[below] {$s$};
			\draw (7.5,1.4) node {p odd};
			\draw (17.5,1.4) node {q odd};
			\draw[line width=1] (5,0) -- (10,0);
			\draw[line width=1] (15,0) -- (20,0);
			\draw[fill= black] (10,0) circle (0.5) node[below] 	{$t$};
			\draw[fill= black] (15,0) circle (0.5)node[below] {$u$};
			\draw[fill= black] (20,0) circle (0.5)node[below] {$v$}; 
			\end{tikzpicture}
	}}
	\begin{proof}
		Lemma \ref{lem:h2 transfer 3 to 2} gives the image of the left hand map. To compute the right hand map we consider the index of spherical subgroups, by Lemma~\ref{lem:transfer on bottom row}. Computing the index of each subgroup as in Lemma~\ref{coset quoients} gives non-zero maps as required.		
	\end{proof}
\end{lem}

\begin{prop}\label{prop:(3,0)entry on E2 for h3}
	Recall from Definition \ref{def:diagrams for h3} the diagrams $\D_{\bullet \bullet }^\square$, $\D{\begin{tikzpicture}[scale=0.08, baseline={(0,0)}]
		\draw[fill= black] (12.5,0) circle (0.6) ;
		\draw (17.5,1.4) node {{\tiny 2r}};
		\draw[line width=1] (15,0) -- (20,0);
		\draw[fill= black] (15,0) circle (0.6);
		\draw[fill= black] (20,0) circle (0.6); 
		\end{tikzpicture}}$ and $\D_{\bf{A}_3}$. Then the $E^2_{3,0}$ of the isotropy spectral sequence for~$(W,S)$ is given by
	\begin{equation*}
	E^2_{3,0} =
	H_1(\D_{\bullet \bullet }^\square;\Z_2) \oplus H_0(\D{\begin{tikzpicture}[scale=0.08, baseline={(0,0)}]
		\draw[fill= black] (12.5,0) circle (0.6) ;
		\draw (17.5,1.4) node {{\tiny 2r}};
		\draw[line width=1] (15,0) -- (20,0);
		\draw[fill= black] (15,0) circle (0.6);
		\draw[fill= black] (20,0) circle (0.6); 
		\end{tikzpicture}};\Z_2) \oplus H_0(\D_{\bf{A}_3};\Z_2) \oplus \big(\underset{\substack{{\scriptscriptstyle W({\bf H}_3)\subseteq W}\\{\scriptscriptstyle W({\bf B}_3)\subseteq W}}}{\oplus}\Z_2 \big).
	\end{equation*}
	\begin{proof}
		Splitting the $d^1$ differentials of Lemma \ref{lem:d1 at (3,0) for h3} as in Proposition \ref{prop:h2 splitting calc}, we can equate the homology of the sequence in Lemma \ref{lem:d1 at (3,0) for h3} to the components on the right-hand side above.
	\end{proof}
\end{prop}

\subsection{Further differentials are zero}\label{sec:h3 further differentials zero}

Recall the isotropy spectral sequence for the Davis complex associated to a Coxeter system $(W,S)$, given in Figure \ref{fig:isotropy ss}. Then on the $p+q=3$ diagonal the spectral sequence has $E^2$ page as shown in Figure~\ref{fig:isotropy ss E2 page}.
	\begin{figure}
		\begin{tikzpicture}
		\matrix (m) [matrix of math nodes,
		nodes in empty cells,nodes={minimum width=5ex,
			minimum height=5ex,outer sep=-5pt},
		column sep=1ex,row sep=1ex]{
			3 	&  	\color{teal}{0}	& \cdots\\
			2 	& 	0	& \color{teal}{A}  & ?& \cdots \\
			1   & 0  &  0 & \color{teal}{B} &? & \cdots\\
			0   &  	\Z	& ? &  ? & \color{teal}{C}  & ?&\\
			\quad\strut &   0  			&  1  &  2  &  3 & 4  \\};
		\draw[thick] (m-1-1.east) -- (m-5-1.east) ;
		\draw[thick] (m-5-1.north) -- (m-5-6.north) ;
		\end{tikzpicture}
		\begin{eqnarray*}
		A&=&H_0(\D_{\operatorname{odd}};\Z_2)\\
		B&=&H_0(\D_{\bullet \bullet };\Z_2) \oplus H_0(\D_{{\bf A}_2};\Z_3)\oplus\big(\underset{m(s,t)>3,\neq \infty}{\oplus}\Z_{m(s,t)} \big)\\
	    C&=&H_1(\D_{\bullet \bullet }^\square;\Z_2)  \oplus H_0(\D{\begin{tikzpicture}[scale=0.08, baseline={(0,0)}]
	    	\draw[fill= black] (12.5,0) circle (0.6) ;
	    	\draw (17.5,1.4) node {{\tiny 2r}};
	    	\draw[line width=1] (15,0) -- (20,0);
	    	\draw[fill= black] (15,0) circle (0.6);
	    	\draw[fill= black] (20,0) circle (0.6); 
	    	\end{tikzpicture}};\Z_2) \oplus H_0(\D_{\bf{A}_3};\Z_2) \oplus \big(\underset{\substack{{\scriptscriptstyle W({\bf H}_3)\subseteq W}\\{\scriptscriptstyle W({\bf B}_3)\subseteq W}}}{\oplus}\Z_2 \big).
		\end{eqnarray*}
		\caption{The $E^2$ page of the isotropy spectral sequence for the Davis complex of a Coxeter system $(W,S)$.}
		\label{fig:isotropy ss E2 page}
	\end{figure}

 The $E^{\infty}$ page of this spectral sequence gives us filtration quotients for $H_3(W;\Z)$ on this diagonal. The arguments in this section shows that all possible further differentials to and from this diagonal are zero. Since the spectral sequence is first quadrant from Figure \ref{fig:isotropy ss E2 page} there are only 3 possible further differentials that may affect the $p+q=3$ diagonal:
\begin{center}
	 (1) $d^2: E^2_{3,1}\to A$ \hspace{1cm}
	 (2) $d^2: E^2_{4,0}\to B$\hspace{1cm}
	 (3) $d^3: E^3_{4,0}\to E^3_{1,2}$.
\end{center}
We first prove two lemmas which will reduce the cases for which we compute differentials originating at $E^r_{4,0}$ in cases~(2) and~(3). Let~$W_A$ and~$W_B$ be non-trivial finite groups, and the size of their generating sets $S_A$ and $S_B$ sum to 4. Denote the isotropy spectral sequence for~$W_A \times W_B$ by~$E(A\times B)$. Then the $E^1_{4,0}$ term in the spectral sequence is 
$$
E^1_{4,0}=H_0(W_A\times W_B;\Z_{A \sqcup B}).
$$
\begin{lem}\label{lem:kunneth1}
	With notation as above, the possible $d^2$ and $d^3$ differentials originating at $E^r_{4,0}$, for $r=2$ or $r=3$, in the spectral sequence $E(A \times B)$ are zero.
	\begin{proof}
		By the K\"{u}nneth theorem for group homology (see e.g.~\cite{Brown}) we have the short exact sequence:		
		{\begin{multline*}
			{0 \to \bigoplus_{i+j=k} H_i(W_A; \Z_A) \underset{\Z}{\otimes} H_j(W_B;\Z_B) \overset{\times}{\to} H_k(W_A \times W_B; \Z_{A \sqcup B})}\\{ \to \bigoplus_{i+j=k-1} \Tor_1^{\Z}(H_i(W_A;\Z_A), H_j(W_B;\Z_B)) \to 0}
			\end{multline*}		}
		
		\noindent since $\Z_{A} \otimes \Z_{B} \cong \Z_{A \sqcup B}$. When $k=0$ the torsion term is zero,
		hence
		$$
		H_0(W_A; \Z_A) \otimes_{\Z} H_0(W_B;\Z_B) \overset{\cong}{\to} H_0(W_A \times W_B; \Z_{A\sqcup B}).
		$$
		\noindent By Theorem \ref{lem:pairing May} there is a pairing 
		$$\Phi_*:E(A)\otimes E(B) \to E(A\times B)
		$$  
		which is given on individual summands of the~$E^1$ terms by the K\"{u}nneth map. Since~$E^1_{4,0}(A\times B)$ has only one summand, $\Phi_*$ is given by the K\"{u}nneth map above, which is an isomorphism. Let $\mid S_A \mid=\alpha$ and $\mid S_B \mid =\beta$ and recall $\alpha + \beta=4$. Then under the pairing $\Phi_*$ all cycles in $E^1_{4,0}(A\times B)$ correspond to a pair of cycles:		
		$$ E^1_{\alpha,0}(A)\otimes E^1_{\beta, 0}(B)\overset{\cong}{\to} E^1_{4,0}(A\times B).$$
		
		It follows that all $d^1$ differentials from $E^1_{4,0}(A\times B)$ are described via the Leibniz rule by differentials from $E^1_{\alpha,0}(A)$ and  $E^1_{\beta,0}(B)$. Therefore the kernel of $d^1$ from~$E^1_{4,0}(A\times B)$ is given by a pairing of elements in the kernel of $d^1$ from $E^1_{\alpha,0}(A)$ and the kernel of $d^1$ from $E^1_{\beta,0}(B)$, and so the K\"{u}nneth map is onto on the $E^2$ page:		 
		$$ E^2_{\alpha,0}(A)\otimes E^2_{\beta, 0}(B)\rightarrow E^2_{4,0}(A\times B)$$		
		and the $d^2$ differentials from $E^2_{4,0}(A\times B)$ are again defined via the Leibniz rule. Since $\alpha$ and $\beta$ are both less than 4, the $d^2$ differentials in~$E(A)$ and~$E(B)$ arise at~$E^2_{p,0}$ where $p<4$. But all possible targets of a $d^2$ differential from such an $E^2_{p,0}$ are zero (consider Figure \ref{fig:isotropy ss E2 page}). Thus the further differentials mapping from $E^2_{4,0}(A\times B)$ are zero.
		
		The $d^2$ differential with target $E^2_{4,0}(A\times B)$ originates at a $0$ group, since the spectral sequence is first quadrant. Since the $d^2$ with source~$E^2_{4,0}(A\times B)$  is also zero,~$E^2_{4,0}(A\times B)=E^3_{4,0}(A\times B)$. By a similar argument~$E^2_{\alpha,0}(A)=E^3_{\alpha,0}(A)$ and~$E^2_{\beta, 0}(B)=E^3_{\beta, 0}(B)$. It follows that the K\"{u}nneth map is also onto on the~$E^3$ page and therefore by the same argument as the~$d^2$ case, the $d^3$ differential originating at $E^3_{4,0}(A\times B)$ is zero.
	\end{proof}
\end{lem}

\begin{lem}\label{lem:kunneth}
	Consider a differential $d^2$ or $d^3$ originating from a summand in~$E^r_{4,0}$ for $r=2$ or $r=3$. If the corresponding cycle at the $E^1_{4,0}$ term occurs in a summand~$H_0(W_A\times W_B;\Z_{A\sqcup B})$, for $W_A$ and $W_B$ non-trivial subgroups of $W$, then the $d^2$ or $d^3$ differential is zero.
	\begin{proof}
		By Lemma \ref{lem:inculsion of gps gives inclusion of ss coxeter}, the inclusion of groups $W_A \times W_B\hookrightarrow W$ gives an inclusion of spectral sequences on the $E^1$ page
		$
		E^1(A \times B)\hookrightarrow E^1(W).
		$
		Therefore differentials mapping from cycles corresponding to the $H_0(W_A\times W_B;\Z_{A\sqcup B})$ summand at position $E^1_{4,0}$ in $E(W)$ will be induced via this map by differentials in $E(A \times B)$. From Lemma \ref{lem:kunneth1} the $d^2$ and $d^3$ differentials originating at the $E^r_{4,0}$ position are zero in~$E(A\times B)$.
	\end{proof}
\end{lem}


We therefore only need to consider differentials originating at the $E^r_{4,0}$ components for~$r=2$ or~$r=3$, which correspond to $H_0(W_T; \Z_T)$ summands of $E^1_{4,0}$ for~$W_T$ irreducible groups, namely for $W_T$ of type ${\bf A}_4, {\bf B}_4, {\bf D}_4, {\bf F}_4$ and ${\bf H}_4$. As in the previous sections we denote the generator of $H_0(W_T;\Z_T)=\Z_2$ by $\D_{W_T}$.

\begin{lem}\label{lem:maps at E4,0}
	The $d^1$ differentials on the $E^1$ page at the $E^1_{4,0}$ position for the summands $H_0(W_T;\Z_T)$ corresponding to Coxeter groups of type ${\bf A}_4, {\bf B}_4, {\bf D}_4, {\bf F}_4$ and~${\bf H}_4$ are non-zero in the single case shown below
	
	\centerline{	\xymatrix@R=-2mm {
			\underset{| T | =3}{\underset{T \in \mathcal{S}}{\oplus}}  H_0(W_T; \Z_T) & \ar[l]_{d^1} \underset{| T | =4}{\underset{T \in \mathcal{S}}{\oplus}}  H_0(W_T; \Z_T)&  \ar[l]_{d^1} \underset{| T | =5}{\underset{T \in \mathcal{S}}{\oplus}}  H_0(W_T; \Z_T) \\
			\begin{tikzpicture}[scale=0.15, baseline=0]
			\draw[fill= black] (5,0) circle (0.5) node[below] {$s$};
			\draw[line width=1] (5,0) -- (15,0);
			\draw[fill= black] (10,0) circle (0.5) node[below] {$t$};
			\draw[fill= black] (15,0) circle (0.5)node[below] {$u$}; 
			\end{tikzpicture} + 	\begin{tikzpicture}[scale=0.15, baseline=0]
			\draw[line width=1] (10,0) -- (20,0);
			\draw[fill= black] (10,0) circle (0.5) node[below] {$t$};
			\draw[fill= black] (15,0) circle (0.5)node[below] {$u$};
			\draw[fill= black] (20,0) circle (0.5)node[below] {$v$}; 
			\end{tikzpicture}& \ar@{|->}[l]	\begin{tikzpicture}[scale=0.15, baseline=0]
			\draw[fill= black] (5,0) circle (0.5) node[below] {$s$};
			\draw[line width=1] (5,0) -- (20,0);
			\draw[fill= black] (10,0) circle (0.5) node[below] {$t$};
			\draw[fill= black] (15,0) circle (0.5)node[below] {$u$};
			\draw[fill= black] (20,0) circle (0.5)node[below] {$v$}; 
			\end{tikzpicture}
	}}
	\begin{proof}
		From Lemma \ref{lem:d1 at (3,0) for h3} we have the maps from the central groups to the left. The finite Coxeter groups with 5 generators for which the ${\bf A}_4, {\bf B}_4, {\bf D}_4, {\bf F}_4$ and ${\bf H}_4$ diagrams are subdiagrams are the groups of type ${\bf A}_5, {\bf B}_5, {\bf D}_5$ and the groups created by taking the product with~${\bf A}_1$. Recall from Lemma \ref{lem:transfer on bottom row} that in this case~$d^1$ is determined by the index of the subgroup. In the case of the product groups, the index of the 4-generator subgroup is $2$ and hence the transfer map is zero. The remaining computations we compute using Python and \cite{Geck}, though formulas for each group size can be found in \cite{Humphreys}. In each case the index of the subgroup is even, hence the transfer map is zero.
	\end{proof}
\end{lem}

\begin{prop}\label{prop:higher d's 0 if trasfers identically to 0}
	If $d^1$ applied to a generator of a summand~$H_q(W_T;\Z_T)$ on the~$E^1$ page is identically zero on the chain level, then the higher differentials which originate at cycles corresponding to this generator on the $E^r$ page are also zero.
	\begin{proof}
		The $d^1$ differential of the isotropy spectral sequence is given by the transfer map on the chain level by Proposition \ref{prop:transfer maps for davis isotropy}. In general higher differentials of the spectral sequence for a double complex are induced by combinations of the differentials on the chain level, and lifting on the chain level. Therefore if the $d^1$ differential is zero on the chain level for the cycle representing a term $E^r_{p,q}$, then the higher differentials will also be zero.
	\end{proof}
\end{prop}

\begin{cor}\label{cor:d2 and d3 B}
	The $d^2$ and $d^3$ differentials originating at $E^r_{0,4}$ for~$r=2$ or~$r=3$ corresponding to cycles on the $E^1_{4,0}$ summands for groups of type ${\bf B}_4, {\bf D}_4, {\bf F}_4$ and~${\bf H}_4$ are zero.
	\begin{proof}
		This is a consequence of Lemma \ref{lem:maps at E4,0}, and Proposition \ref{prop:higher d's 0 if trasfers identically to 0}, if the transfer maps from Lemma \ref{lem:transfer on bottom row} originating at $H_0(W_T;\Z_T)$ for these groups are identically zero on the chain level (and not just zero modulo 2). This is satisfied if, alongside there being an even number of cosets, there are identical numbers of cosets with odd and even length. We use Python \cite{Geck} and compute that there are equal numbers of coset representatives of even and odd length for every three generator subgroup of ${\bf B}_4, {\bf D}_4, {\bf F}_4$ and ${\bf H}_4$.
	\end{proof}
\end{cor}

\noindent The remaining potentially non-zero differentials originating at the $E^r_{0,4}$ position for~$r=2$ or~$r=3$, correspond to cycles on the~$E^1_{4,0}$ summand~$H_0(W({\bf A}_4);\Z_T)$.

\begin{lem} \label{lem: ss for A4}
	The potential~$d^2$ and~$d^3$ differentials originating at the~$E^r_{0,4}$ position for~$r=2$ or~$r=3$ and corresponding to cycles on the~$E^1_{4,0}$ summand~$H_0(W({\bf A}_4);\Z_T)$ are zero.
	\begin{proof}
		If the further differentials were non zero then they would also be non zero in the spectral sequence for~$W({\bf A}_4)$ by Lemma~\ref{lem:inculsion of gps gives inclusion of ss coxeter}. The~$E^2$ page for the Coxeter group~$W({\bf A}_4)$ is given by Figure~\ref{fig:isotropy ss E2 page} with
		\begin{eqnarray*}
		A=0 \hspace{1cm}
		B = \Z_2\oplus \Z_3 \hspace{1cm}
		C= \Z_2.
		\end{eqnarray*}
		\noindent The computation of this is given in \cite[Appendix B]{Boydthesis}.
		The third integral homology of the symmetric group~$S_5$, which is isomorphic to~$W({\bf A}_4)$, is 
		$$
		H_3(W({\bf A}_4);\Z)=\Z_{12}\oplus \Z_2\cong \Z_3 \oplus \Z_4 \oplus \Z_2
		$$
		\noindent which is precisely given if the groups on the $p+q=3$ diagonal of the $E^2$ page are the $E^\infty$ terms, or \emph{filtration quotients} for $H_3(W({\bf A}_4);\Z)$ (there is a non-trivial extension of $\Z_2$ by $\Z_2$ which we will discuss in the following section). Therefore no higher differentials in or out of this diagonal can be non-zero.
	\end{proof}
\end{lem}

\begin{prop}
	The possible $d^2$ and $d^3$ differentials originating at the $E^r_{4,0}$ position for~$r=2$ or~$r=3$ in the spectral sequence are zero.
	\begin{proof}
		This is direct result of Lemma \ref{lem:kunneth}, Corollary \ref{cor:d2 and d3 B} and Lemma \ref{lem: ss for A4}.
	\end{proof}
\end{prop}

\begin{lem}\label{lem:kunneth with tor 1}
	Let $W_T$ and $W_V$ be non-trivial finite Coxeter groups, and the size of their generating sets sum to 3. Then the potential $d^2$ differential originating at the~$E^2_{3,1}$ position is zero.
\begin{proof}
	The group $W_T\times W_V$ must be $W({\bf I}_2(p))\times W({\bf A}_1)$ for $p\geq 2$, by the classification of finite Coxeter groups. 
	
	When $p$ is even, the $E^2$ page for the Coxeter group $W({\bf I}_2(p))\times W({\bf A}_1)$ is given by Figure \ref{fig:isotropy ss E2 page} with
		\begin{eqnarray*}
		A&=&\Z_2\oplus \Z_2\oplus \Z_2 \hspace{1cm}
		B = \Z_2 \oplus \Z_2\oplus \Z_{p} \hspace{1cm}
		C= \Z_2
	\end{eqnarray*}		
	\noindent which is computed in \cite[Appendix B]{Boydthesis}. The third integral homology can be computed via the K\"{u}nneth formula for groups to be
	\begin{eqnarray*}
		H_3(W({\bf I}_2(p))\times W({\bf A}_1); \Z)
		&=& \Z_2 \oplus \Z_2 \oplus  \Z_2\oplus \Z_2 \oplus \Z_{p} \oplus \Z_2 \oplus \Z_2.
	\end{eqnarray*}		
	
	Similarly, when $p$ is odd, the $E^2$ page is given by Figure \ref{fig:isotropy ss E2 page} with
	\begin{eqnarray*}
		A&=&\Z_2\oplus \Z_2 \hspace{1cm}
		B =  \Z_2\oplus \Z_{p} \hspace{1cm}
		C= 0
	\end{eqnarray*}
	
	\noindent and the K\"{u}nneth formula gives the homology to be
	\begin{eqnarray*}
		H_3(W({\bf I}_2(p))\times W({\bf A}_1); \Z)
		&=& \Z_2 \oplus \Z_{2p} \oplus \Z_2.
	\end{eqnarray*}	
	
	\noindent In both cases, the group homology calculated via K\"{u}nneth is precisely given if the groups on the $p+q=3$ diagonal of the $E^2$ page are the $E^\infty$ terms. Therefore no higher differentials in or out of this diagonal can be non-zero.
\end{proof}
\end{lem}

\begin{lem}\label{lem:kunneth with tor for d2 in h3 calc}
	Suppose a $d^2$ differential in the isotropy spectral sequence for $W$ originates at a cycle in $E^2_{3,1}$ represented by a homology class in $E^1_{3,1}$ of a subgroup~$W_T\times W_V$ of $W$ such that neither $W_T$ or $W_V$ is the trivial group. Then this $d^2$ differential is zero.
	\begin{proof}
		This proof mimics Lemma~\ref{lem:kunneth}, using Lemma~\ref{lem:inculsion of gps gives inclusion of ss coxeter}, and Lemma \ref{lem:kunneth with tor 1}.
	\end{proof}
\end{lem}

\begin{prop}
	The possible $d^2$ differential originating at the $E^2_{3,1}$ group in the spectral sequence is zero.
	\begin{proof}
		The $E^2_{3,1}$ entry is calculated by computing the homology of the sequence 
		\centerline{	\xymatrix@R=0.5mm {
							\underset{| T | =2}{\underset{T \in \mathcal{S}}{\oplus}}  H_1(W_T; \Z_T) & \ar[l]_{d^1} \color{teal}{\underset{| T | =3}{\underset{T \in \mathcal{S}}{\oplus}}  H_1(W_T; \Z_T)}&  \ar[l]_{d^1} \underset{| T | =4}{\underset{T \in \mathcal{S}}{\oplus}}  H_1(W_T; \Z_T). }}
		Recall the left hand map from Proposition \ref{prop:d1 E^2 for H_1s}.
		The possible $d^2$ differential acts on cycles in summands of the form $H_1(W_T;\Z_T)$ for $|T|=3$.
		
		If~$d^2$ acts on a cycle in the summand $H_1(W({\bf A}_3);\Z_T)=\Z_3$ (from Lemma \ref{prop:h1 of 3 gen groups}), it must map to zero, since the target~$E^2_{1,2}=H_0(\D_{\operatorname{odd}};\Z_2)$ is all 2-torsion. 
		
		If~$d^2$ acts on a cycle in the summand $H_1(W_T;\Z_T)$ for~$W_T$ the group~$W({\bf B}_3)$ or~$W({\bf H}_3)$ it will map to zero, as the representing cycles transfer identically to zero on the chain level by the proof of Lemma \ref{prop:d1 E^2 for H_1s}, so we apply Proposition \ref{prop:higher d's 0 if trasfers identically to 0}. 		
		
		Lemma \ref{lem:kunneth with tor for d2 in h3 calc} covers the final cases where the $d^2$ acts on a cycle in the summand $H_1(W_T;\Z_T)$ for $W_T=W({\bf I}_2(p))\times W({\bf A}_1)$ for $ p\geq 2$.
	\end{proof}
\end{prop}

\subsection{Extension problems}\label{sec:extension problem}

Since all further differentials at the $p+q=3$ diagonal are zero, the $E^2$ page shown in Figure \ref{fig:isotropy ss E2 page} gives the limiting, or $E^\infty$, terms on this diagonal. The spectral sequence on this diagonal converges to filtration quotients of~$H_3(W;\Z)$, so we consider possible extensions on this diagonal. That is there is a filtration of $H_3(W;\Z)$
$$
F_0\subseteq F_1 \subseteq F_2 \subseteq F_3 = H_3(W; \Z)
$$
where 
$E^\infty_{0,3}=F_0$, $E^\infty_{1,2}=F_1 / F_0$, $E^\infty_{2,1}=F_2 / F_1$ and $E^\infty_{3,0}=F_3 / F_2$. We have $F_0=0$ and so $E^\infty_{1,2}=F_1$. 

\begin{prop}
	The group $F_1=A=H_0(\D_{\operatorname{odd}};\Z_2)$ splits off as a direct summand of $H_3(W;\Z)$.
	\begin{proof}
		Consider a homomorphism $\psi$ from a Coxeter group $W$ with generating set~$S$ to the cyclic subgroup of order two generated by $t$ in $S$, which we denote~$W_t$. If $s_1$ and $s_2$ in~$S$, satisfy $m(s_1,s_2)$ is odd we require $\psi(s_1)=\psi(s_2)$, whereas if $m(s_1,s_2)$ is even there is no requirement on $\psi$. A summand of $$A=F_1=H_0(\D_{\operatorname{odd}};\Z_2)=\underset{\pi_0(\D_{\operatorname{odd}})}{\bigoplus} \Z_2$$ \noindent is represented by a vertex of $\D_W$. For the vertex $t$, denote the corresponding summand of $A$ by $\Z_2(t)$. We define the homomorphism $\psi$ from $W$ to $W_t$ to be zero on all but one of the connected components of $\D_{\operatorname{odd}}$, namely the $t$ component.
		\ms
		\begin{eqnarray*}
			\psi: W &\to& W_t\\
			s &\mapsto& \begin{cases}
				t & \text{if } s \text{ and } t \text{ are in the same component of }\pi_0(\D_{\operatorname{odd}})\\
				e & \text{otherwise.}
			\end{cases}
		\end{eqnarray*} 
		\ms
		
		\noindent Then the map $\psi$ induces a map $\psi_*$ which fits into the following diagram  
		\ms
		
		\centerline{\xymatrix{					
				\Z_2(t) \ar@{^{(}->}[r]\ar[rrrd]_{\id}&A \ar@{^{(}->}[r]&H_3(W;\Z)\ar[r]^-{\psi_*}& H_3(W_t;\Z)\ar@{=}[d]\\
				&&&\Z_2
		}}
		\ms
		
		\noindent where $H_3(W_t;\Z)=\Z_2$ is computed by noting that the $E^{\infty}$ page of the isotropy spectral sequence for $W_t$ has only the group $H_0(\D_{\operatorname{odd}};\Z_2)=\Z_2(t)$ on the $p+q=3$ diagonal.
		The inclusion map $A \hookrightarrow H_3(W;\Z)$ comes from the fact that $A=F_1$ and so is a subgroup of~$H_3(W;\Z)$. The identity isomorphism gives us that $H_3(W;\Z)$ splits as 
		$$
		H_3(W;\Z)=\Z_2(t) \oplus \ker(\psi_*)
		$$
		\noindent and so there are no non-trivial extensions involving the $\Z_2(t)$ summand of $A$. Repeating this argument over all summands gives that there are no non-trivial extensions involving $A$ and so $A=F_1$ splits off in $H_3(W;\Z)$, as required.
	\end{proof}
\end{prop}

We therefore have the filtration 
$$
0\subseteq F_1 \subseteq F_2 \subseteq F_3 = H_3(W; \Z)=F_1\oplus F_3'
$$
and we let $F_2=F_1\oplus F_2'$ and $F_3=F_1\oplus F_3'$. It follows that $E^\infty_{2,1}=B=F_2 / F_1=F_2'$ and $E^\infty_{0,3}=C=F_3/F_2=F_3'/F_2'$, so $F_3'$ fits into the following exact sequence

\centerline{\xymatrix{
		0 \ar[r]&F_2' \ar@{=}[d]\ar[r]& F_3' \ar@{=}[d]\ar[r]& F_3'/F_2' \ar@{=}[d]\ar[r]& 0\\
		0 \ar[r]& B \ar[r]& F_3' \ar[r]& C \ar[r]& 0
}}
\noindent i.e.~ $F_3'$ is an extension of $C$ by $B$.

\begin{lem}\label{lem:ex1}
	There exist no non-trivial extensions between the $H_0(\D{\begin{tikzpicture}[scale=0.08, baseline={(0,0)}]
		\draw[fill= black] (12.5,0) circle (0.6) ;
		\draw (17.5,1.4) node {{\tiny 2r}};
		\draw[line width=1] (15,0) -- (20,0);
		\draw[fill= black] (15,0) circle (0.6);
		\draw[fill= black] (20,0) circle (0.6); 
		\end{tikzpicture}};\Z_2)$ summand of C and the groups at $B$ in the spectral sequence of Figure \ref{fig:isotropy ss E2 page}.
	\begin{proof}
		A summand of $H_0(\D{\begin{tikzpicture}[scale=0.08, baseline={(0,0)}]
			\draw[fill= black] (12.5,0) circle (0.6) ;
			\draw (17.5,1.4) node {{\tiny 2r}};
			\draw[line width=1] (15,0) -- (20,0);
			\draw[fill= black] (15,0) circle (0.6);
			\draw[fill= black] (20,0) circle (0.6); 
			\end{tikzpicture}};\Z_2)$ is represented by a vertex in $\D{\begin{tikzpicture}[scale=0.08, baseline={(0,0)}]
			\draw[fill= black] (12.5,0) circle (0.6) ;
			\draw (17.5,1.4) node {{\tiny 2r}};
			\draw[line width=1] (15,0) -- (20,0);
			\draw[fill= black] (15,0) circle (0.6);
			\draw[fill= black] (20,0) circle (0.6); 
			\end{tikzpicture}}$ corresponds to an ${\bf I}_2(2p)\sqcup {\bf A}_1$ ($p\geq1$) subdiagram present in $\D_W$. We compute the spectral sequence for the Coxeter group $V=W({\bf I}_2(2p))\times W({\bf A}_1)$ corresponding to this diagram, and note that by Lemma \ref{lem:inculsion of gps gives inclusion of ss coxeter} the inclusion of the subgroup $V$ into the group $W$ induces a map of spectral sequences. Therefore if there is a trivial extension in the spectral sequence for $V$ corresponding to the ${\bf I}_2(2p)\sqcup {\bf A}_1$ summand of $H_0(\D{\begin{tikzpicture}[scale=0.08, baseline={(0,0)}]
			\draw[fill= black] (12.5,0) circle (0.6) ;
			\draw (17.5,1.4) node {{\tiny 2r}};
			\draw[line width=1] (15,0) -- (20,0);
			\draw[fill= black] (15,0) circle (0.6);
			\draw[fill= black] (20,0) circle (0.6); 
			\end{tikzpicture}};\Z_2)$, this extension will be trivial in the spectral sequence for $W$. This is because the splitting of the extension sequence in $E(V)$ will give a splitting of the extension sequence in $E(W)$, under the map of spectral sequences.
		 The $E^\infty$ page for the Coxeter group $V$ is given by Figure \ref{fig:isotropy ss E2 page} with
		\begin{eqnarray*}
			A&=&\Z_2\oplus \Z_2\oplus \Z_2 \hspace{1cm}
			B = \Z_2 \oplus \Z_2\oplus \Z_{2p} \hspace{1cm}
			C= \Z_2
		\end{eqnarray*}		
		\noindent which is computed in \cite[Appendix B]{Boydthesis}. Therefore
		$$H_3(V;\Z)=F_3'\oplus F_1=F_3'\oplus (\Z_2\oplus \Z_2 \oplus \Z_2)$$ where $F_3'$ is an extension of $\Z_2$ by $\Z_2 \oplus\Z_2 \oplus  \Z_{2p}$.
		
		The third integral homology of $V$ can be computed via the K\"{u}nneth formula for groups to be
		\begin{eqnarray*}
			H_3(W({\bf I}_2(2p))\times W({\bf A}_1); \Z)
			&=& \Z_2 \oplus \Z_2 \oplus  \Z_2\oplus \Z_2 \oplus \Z_{2p} \oplus \Z_2 \oplus \Z_2.
		\end{eqnarray*}		
		
		\noindent Therefore we see that $F_3'=\Z_2 \oplus \Z_2 \oplus \Z_2 \oplus \Z_{2p}$ and it follows that there is no non-trivial extension between the $H_0(\D{\begin{tikzpicture}[scale=0.08, baseline={(0,0)}]
			\draw[fill= black] (12.5,0) circle (0.6) ;
			\draw (17.5,1.4) node {{\tiny 2r}};
			\draw[line width=1] (15,0) -- (20,0);
			\draw[fill= black] (15,0) circle (0.6);
			\draw[fill= black] (20,0) circle (0.6); 
			\end{tikzpicture}};\Z_2)$ component of $C$ and $B$.
	\end{proof}
\end{lem}

\begin{lem}\label{lem:ex2}
	The extension between the $H_0(\D_{\bf{A}_3};\Z_2)$ summand in C and the  $H_0(\D_{\bullet \bullet };\Z_2)$ summand in $B$ is non-trivial.
	\begin{proof}
		A summand of $H_0(\D_{\bf{A}_3};\Z_2)$ is represented by a vertex of $\D_{\bf{A}_3}$, corresponding to an ${\bf A}_3$ subdiagram present in $\D_W$. The~$E^\infty$ page of spectral sequence for the subgroup $V=W({\bf A}_3)$ corresponding to this diagram is given by Figure \ref{fig:isotropy ss E2 page} with
		\begin{eqnarray*}
			A&=&\Z_2\hspace{1cm}
			B = \Z_2  \oplus\Z_{3}\hspace{1cm}
			C= \Z_2
		\end{eqnarray*}
		\noindent which is computed in \cite[Appendix B]{Boydthesis}.
		Therefore
		$$H_3(V;\Z)=F_3'\oplus F_1=F_3'\oplus \Z_2$$ where $F_3'$ is an extension of $\Z_2$ by $\Z_2 \oplus \Z_3$. Recall that $V$ is isomorphic to the symmetric group $S_4$, and~$H_3(S_4;\Z)=\Z_{12} \oplus \Z_2$. The unique extension which will obtain this result is the following:		
		$$
		0 \to \Z_2 \oplus \Z_3 \to \Z_4 \oplus \Z_3 \to \Z_2 \to 0
		$$		
		\noindent giving $H_3(V;\Z)=\Z_4\oplus \Z_3\oplus \Z_2=\Z_{12} \oplus \Z_2$. By Lemma \ref{lem:inculsion of gps gives inclusion of ss coxeter} the inclusion of subgroup $V$ into group $W$ gives a map of spectral sequences, under which the extension sequence above is mapped as follows:
		
		\centerline{\xymatrix{
				0 \ar[r] &\ar@{^{(}->}[d] \Z_2 \oplus \Z_3 \ar[r] &\Z_4 \oplus \Z_3 \ar@{.>}[d] \ar[r] &\Z_2 \ar@{^{(}->}[d]\ar[r] &0 \\		
				0 \ar[r]& B \ar[r]& F_3' \ar[r]& C \ar[r] &0.
		}}
		
		\noindent Therefore the extension in $E(V)$ corresponding to the ${\bf A}_3$ summand of $H_0(\D_{\bf{A}_3};\Z_2)$ is present in the spectral sequence for $W$. It follows that there exists a non-trivial extension from each summand of $H_0(\D_{\bf{A}_3};\Z_2)$ to $H_0(\D_{\bullet \bullet };\Z_2)$.
	\end{proof}
\end{lem}
\begin{defn}\label{def:extension matrix}
	For a Coxeter system $(W,S)$, let $I=\pi_0(\D_{\bullet \bullet })$, $J=\pi_0(\D_{\bf{A}_3})$, let the connected component of a vertex $\{s,u\}$ in $\pi_0(\D_{\bullet \bullet })$ be denoted $[\{s,u\}]$ and the connected component of a vertex $\{s,t,u\}$ in $\pi_0(\D_{\bf{A}_3})$ be denoted $[\{s,t,u\}]$. We define the \emph{extension matrix} $X_W$ to be the $|I|$ by $|J|$ matrix with entries
	
	$$X(i,j)=\begin{cases}1 & \text{ if }i=[\{s,u\}] \text{ and } j=[\{s,t,u\}] \\ 0 & \text{ otherwise}.
	\end{cases}
	$$
\end{defn}

\begin{lem}\label{lem:extension is given by extension matrix}
	The extension of $H_0(\D_{\bf{A}_3};\Z_2)$ by $H_0(\D_{\bullet \bullet };\Z_2)$ in the spectral sequence is completely determined by the extension matrix $X_W$ defined in Definition \ref{def:extension matrix}. The extension sequence in question is
	
	\centerline{\xymatrix@R=3mm{
			0 \ar[r]&H_0(\D_{\bt\bt};\Z_2) \ar@{=}[d] \ar[r]&\ar@{=}[d] Y \ar[r]&\ar@{=}[d]H_0(\D_{\bf{A}_3};\Z_2) \ar[r]& 0\\
			0 \ar[r]&\underset{\pi_0(\D_{\bullet \bullet })}{\bigoplus} \Z_2 \ar[r]& 
			Y \ar[r]&\underset{\pi_0(\D_{\bf{A}_3})}{\bigoplus} \Z_2 \ar[r]& 0.
	}}

	\noindent and the entry $X(i,j)$ of $X_W$ dictates whether the extension between the $i^{th}$ $\Z_2$ on the left and $j^{th}$ $\Z_2$ on the right is trivial (if $X(i,j)=0$) or $\Z_4$ (if $X(i,j)=1$).
	
	\begin{proof}
		For two finite indexing sets $I$ and $J$, the extensions of $\underset{J}{\oplus} \Z_2$ by $\underset{I}{\oplus} \Z_2$ are classified by \vspace{-.35cm}
		\begin{eqnarray*}
			\operatorname{Ext}(\underset{I}{\oplus} \Z_2,\underset{J}{\oplus} \Z_2)&=&\underset{I}{\oplus}\underset{J}{\oplus}\operatorname{Ext}(\Z_2, \Z_2)
			=\underset{I}{\oplus}\underset{J}{\oplus} \Z_2.
		\end{eqnarray*}
		\noindent Under this classification, an extension is given by an~$I\times J$ matrix $X$ with entries~$X(i,j)$ in~$\Z_2$. The $X(i,j)$ entry is zero if the restriction to these summands in the extension sequence is trivial, and $1$ if the extension is the non-trivial extension of~$\Z_2$ by $\Z_2$ giving $\Z_4$.

		Consider the extension sequence. By Lemma \ref{lem:ex2}, we know that the projection on the right to a $\Z_2$ summand $[\{s,t,u\}]$ in $\pi_0(\D_{\bf{A}_3})$ is the non-trivial extension by the $\Z_2$ summand~$[\{s,u\}]$ in $\pi_0(\D_{\bullet \bullet })$. Let $I=\pi_0(\D_{\bullet \bullet })$ and $J=\pi_0(\D_{\bf{A}_3})$, then the matrix $X$ is precisely $X_W$ from Definition \ref{def:extension matrix}.
	\end{proof}	
\end{lem}

\begin{lem}\label{lem:ex3}
	There exist no non-trivial extensions between the $$\oplus \big(\underset{\substack{{\scriptscriptstyle W({\bf H}_3)\subseteq W}\\{\scriptscriptstyle W({\bf B}_3)\subseteq W}}}{\oplus}\Z_2 \big)$$ \noindent summand of $C$ and the groups at $B$ in the spectral sequence of Figure \ref{fig:isotropy ss E2 page}.  
	\begin{proof}
		We recall that subdiagrams of the form ${\bf H}_3$ and ${\bf B}_3$ in $\D_W$ represent these summands of $C$. We compute the spectral sequence for the groups corresponding to these diagrams, and compare to the third homology of the corresponding group~$W({\bf H}_3)$ or $W({\bf B}_3)$ as computed using the De Concini - Salvetti resolution \cite{DCSal}. Through these comparisons we observe that there are no non-trivial extensions present, as in the proof of Lemma \ref{lem:ex1}. These calculations are in \cite[Appendix B]{Boydthesis}.
\end{proof}		
\end{lem}

\begin{lem}\label{lem:extension problem, class for D squares}
	A class $H_1(\D^\square_{\bt \bt}; \Z_2)$ in $C$ exists only when the spectral sequence is calculated for a Coxeter system $(W,S)$ for which $\D_W$ has a subdiagram of the form $Y \sqcup {\bf A}_1$ where $Y$ is a 1-cycle in the Coxeter diagram $\D_{\operatorname{odd}}$. That is a class in~$H_1(\D^\square_{\bt \bt}; \Z_2)$ is represented in $\D_W$ by a loop containing only odd edges, along with a vertex disjoint from this loop.
	\begin{proof}
		Suppose vertices~$\{t_1,\ldots, t_k\}$ of~$\D_W$ represent a 1-cycle in~$D_{\operatorname{odd}}$ and the  vertex~$s$ is disjoint. Then~$\{(t_1,s),\ldots , (t_k,s)\}$ represents a 1-cycle in~$\D^\square_{\bt \bt}$. To show that all classes in~$H_1(\D^\square_{\bt \bt}; \Z_2)$ are represented by cycles of this form, suppose that~$\{(x_1,y_1), \ldots , (x_p,y_p)\}$ represents a 1-cycle in~$\D^\square_{\bt \bt}$. Without loss of generality, suppose~$x_1=x_2$. Since there exists an edge between~$(x_1,y_1)$ and~$(x_1,y_2)$ in~$\D_{\bt \bt}$,~$m(y_1,y_2)$ must be odd. Now either~$x_1=x_3$ or~$y_2=y_3$; suppose~$y_2=y_3$. It follows that~$m(x_1,x_3)$ is odd, so in~$\D_W$ there is a subdiagram of the form
		\begin{center}
			\begin{tikzpicture}[scale=0.2, baseline=0]
			\draw[fill= black] (5,0) circle (0.4) node[below] {$x_1$};
			\draw (7.5,1.4) node {odd};
			\draw (17.5,1.4) node {odd};
			\draw[line width=1] (5,0) -- (10,0);
			\draw[line width=1] (15,0) -- (20,0);
			\draw[fill= black] (10,0) circle (0.4) node[below] 	{$x_3$};
			\draw[fill= black] (15,0) circle (0.4)node[below] {$y_1$};
			\draw[fill= black] (20,0) circle (0.4)node[below] {$y_2$}; 
			\end{tikzpicture}.
		\end{center}
		\noindent It follows in the diagram $\D_{\bt \bt}$ there is a subdiagram
		\begin{center}
			\begin{tikzpicture}[scale=0.15, baseline=0]
			\draw[fill= black] (5,0) circle (0.5) node[above] {$(x_1, y_1)$};
			\draw[line width=1] (5,0) -- (15,0);
			\draw[line width=1] (5,-10) -- (15,-10);			
			\draw[line width=1] (5,0) -- (5,-10);			
			\draw[line width=1] (15,0) -- (15,-10);
			\draw[fill= black] (15,0) circle (0.5) node[above] 	{$(x_3,y_1)$};
			\draw[fill= black] (5,-10) circle (0.5)node[below] {$(x_1,y_2)$};
			\draw[fill= black] (15,-10) circle (0.5)node[below] {$(x_3,y_2)$}; 
			\end{tikzpicture}
		\end{center}
		\noindent and since this is a square, it is a 2-cell in $\D^\square_{\bt \bt}$. Therefore in~$H_1(\D^\square_{\bt \bt}; \Z_2)$ the cycle~$\{(x_1,y_1), (x_1, y_2), (x_3,y_2), (x_3,y_1)\}$ is a boundary. It follows that replacing the sub-cycle~$\{(x_1,y_1), (x_1, y_2), (x_3,y_2)\}$ of $\{(x_1,y_1), \ldots , (x_p,y_p)\}$ with the vertex~$\{(x_3,y_1)\}$ gives representatives of the same class in $H_1(\D^\square_{\bt \bt}; \Z_2)$, and the original cycle becomes~$\{(x_3,y_1), (x_4,y_4) \ldots , (x_p,y_p)\}$. Without loss of generality, we can now assume that $x_3=x_4$ and we return to the start of the analysis of the cycle. By reiterating this procedure we build a cycle equivalent, via boundaries, to~$\{(x_1,y_1),\ldots, (x_k,y_k)\}$ and where $x_1=x_i$ for all $i$. This is exactly a subdiagram of the form $Y \sqcup {\bf A}_1$ in the Coxeter diagram $\D_W$, where $Y$ is a loop in $\D_{\operatorname{odd}}$.
	\end{proof}
\end{lem}

\begin{lem}\label{lem:homology of class for d squares extension}
	Let $W=W(Y)\times W({\bf A}_1)$ be a Coxeter group such that~$Y$ represents a 1-cycle in $\D_{\operatorname{odd}}$, then for some~$0<m$ in~$\mathbb{N}$
	$$
	H_3(W;\Z)\cong H_3(W(Y);\Z) \oplus \Z_2^m.
	$$
	\begin{proof} By the K\"{u}nneth formula for group homology,
		\begin{eqnarray*}
			H_3(W;\Z)&\cong& H_3(W(Y);\Z) \oplus \Z_2 \oplus H_2(W(Y);\Z) \oplus H_1(W(Y);\Z)
		\end{eqnarray*}	 
	and since the first and second integral homologies of any Coxeter group are all 2-torsion the result follows.
	\end{proof}
\end{lem}

\begin{prop}\label{prop:no extensions in class of d squares}
	When $W=W(Y)\times W({\bf A}_1)$ such that $Y$ represents a 1-cycle in $\D_{\operatorname{odd}}$, there are no non-trivial extensions between the $H_1(\D^\square_{\bt \bt}; \Z_2)$ component in~$C$ and $B$.	
	\begin{proof}
		We note that should non-trivial extensions exist, the homology $H_3(W;\Z)$ would have at least one more summand with torsion greater than 2-torsion, in comparison to the the homology $H_3(W(Y);\Z)$. This is due to the fact that $H_1(\D^\square_{\bt \bt}; \Z_2)$ is zero in the spectral sequence for $H_3(W(Y);\Z)$ so the extension would not occur here. We also note that transitioning from~$W(Y)$ to~$W$ does not alter any non trivial extensions in the spectral sequence for $W(Y)$ between the summand~$H_0(\D_{\bt \bt};\Z_2)$ and~$H_0(\D_{{\bf A}_3};\Z_2)$. From Lemma \ref{lem:homology of class for d squares extension} we have that $H_3(W;\Z)$ has no summands with higher than 2-torsion that do not also appear in $H_3(W(Y);\Z)$.
	\end{proof}
\end{prop}

\begin{lem}\label{lem:ex4}
	There exist no non-trivial extensions from the $H_1(\D^\square_{\bt \bt}; \Z_2)$ component of $C$ to $B$.
	\begin{proof}
		A class of $H_1(\D^\square_{\bt \bt}; \Z_2)$ is represented by a subgroup with diagram of the form~$\D_W=Y \sqcup {\bf A}_1$ such that $Y$ represents a 1-cycle in $\D_{\operatorname{odd}}$, by Lemma~\ref{lem:extension problem, class for D squares}. By Proposition \ref{prop:no extensions in class of d squares} no non-trivial extensions exist between this class and $B$ in the spectral sequence for the representing subgroup. Therefore by similar argument to Lemma \ref{lem:ex1} there are no non-trivial extensions from this class.
	\end{proof}
\end{lem}

\subsection{Proof of Theorem {\ref{thm:INTRO THM B}}}\label{sec:h3 proof of theorem}

\begin{thm}
	Given a finite rank Coxeter system $(W,S)$ there is an isomorphism
\begin{eqnarray*}
	H_3(W; \Z) &\cong& H_0(\D_{\operatorname{odd}};\Z_2)  \oplus H_0(\D_{{\bf A}_2};\Z_3) \oplus \big(\underset{3<m(s,t)< \infty}{\oplus}\Z_{m(s,t)} \big)\\&&\oplus H_0(\D{\begin{tikzpicture}[scale=0.08, baseline={(0,0)}]
		\draw[fill= black] (12.5,0) circle (0.6) ;
		\draw (17.5,1.4) node {{\tiny 2r}};
		\draw[line width=1] (15,0) -- (20,0);
		\draw[fill= black] (15,0) circle (0.6);
		\draw[fill= black] (20,0) circle (0.6); 
		\end{tikzpicture}};\Z_2) \oplus \big(\underset{\substack{{\scriptscriptstyle W({\bf H}_3)\subseteq W}\\{\scriptscriptstyle W({\bf B}_3)\subseteq W}}}{\oplus}\Z_2 \big)\\&&\oplus \big(H_0(\D_{\bf{A}_3};\Z_2) \bigcirc H_0(\D_{\bullet \bullet };\Z_2)\big)  	\oplus H_1(\D^\square_{\bt \bt};\Z_2)
\end{eqnarray*}
	where each diagram is as in Definition \ref{def:diagrams for h3}, and viewed as a cell complex. In this equation, $\bigcirc$ denotes the non-trivial extension of $H_0(\D_{\bf{A}_3};\Z_2)$ by $H_0(\D_{\bullet \bullet };\Z_2)$ given by the extension matrix $X_W$ in Definition \ref{def:extension matrix}.	
	\begin{proof}	
		The extension problems are solved in Lemmas \ref{lem:ex1}, \ref{lem:ex2}, \ref{lem:ex3} and \ref{lem:ex4}. It follows that the only non-trivial extension is the extension of $H_0(\D_{\bf{A}_3};\Z_2)$ by $H_0(\D_{\bullet \bullet };\Z_2)$, which is determined by the extension matrix $X_W$ of Definition \ref{def:extension matrix} by Lemma \ref{lem:extension is given by extension matrix}.
		
		The computation of the $p+q=3$ diagonal of the isotropy spectral sequence for the Davis complex, alongside the solutions to these extension problems, gives the formula for $H_3(W;\Z)$ as stated in the theorem.
	\end{proof}
\end{thm}

\appendix

\section{Table of results for finite Coxeter groups}\label{appendix:table}

The finite Coxeter groups were classified in the 1930s by Coxeter \cite{Coxeter}. This classification is described in Theorem \ref{prop:classification of finite Coxeter}. We use Theorems \ref{thm:INTRO THM A} and \ref{thm:INTRO THM B} to calculate the second and third integral homology of the finite irreducible Coxeter groups, and give the results in Table \ref{table:homology} below. We include $H_1(W;\Z)$ for completeness.

\begin{table}[h!]
	\begin{tabular}[center]{|c|c|c|c|}
		\hline&&&\\
		$W$ & $H_1(W;\Z)$ & $ H_2(W; \Z) $& $H_3(W,\Z)$ \\&&&\\
		\hline &&&\\[-1em]
		${\begin{array}{c c}{\bf A}_n\\n\geq1\end{array}}$& $\Z_2$ & ${\begin{array}{c r}
			0  &n\leq2\\	 
			\Z_2   &n\geq 3\\\end{array}}$ & ${\begin{array}{c r}
			\Z_2  &n=1\\	 
			\Z_2 \oplus \Z_3  &n=2\\
			\Z_2 \oplus \Z_3 \oplus \Z_4 & n=3,4\\
			\Z_2^2 \oplus \Z_3 \oplus \Z_4& n\geq 5 \\ \end{array}}$\\	\hline&&&\\[-1em]
		${\begin{array}{c c}{\bf B}_n\\n\geq2\end{array}}$& $\Z_2\oplus \Z_2$ & ${\begin{array}{c r}
			\Z_2  &n=2\\	 
			\Z_2\oplus \Z_2   &n= 3\\
			\Z_2\oplus \Z_2 \oplus \Z_2   &n \geq 4\\\end{array}}$ & ${\begin{array}{c r}
			\Z_2^2 \oplus \Z_4 &n=2\\	 
			\Z_2^4 \oplus \Z_3 \oplus \Z_4 &n=3\\
			\Z_2^5 \oplus \Z_3 \oplus \Z_4^2 & n= 4\\
			\Z_2^6 \oplus \Z_3 \oplus \Z_4^2 & n= 5\\
			\Z_2^7 \oplus \Z_3 \oplus \Z_4^2 & n\geq 6\\
			\end{array}}$\\	\hline&&&\\[-1em]
		${\begin{array}{c c}{\bf D}_n\\n\geq4\end{array}}$& $\Z_2$ & ${\begin{array}{c r} 
			\Z_2\oplus \Z_2 \\\end{array}}$ & ${\begin{array}{c r}
			\Z_2^2 \oplus \Z_3 \oplus \Z_4^3 & n= 4\\				
			\Z_2^2 \oplus \Z_3 \oplus \Z_4^2 & n= 5\\		
			\Z_2^3 \oplus \Z_3 \oplus \Z_4^2 & n\geq 6\\
			\end{array}}$\\	\hline&&&\\[-1em]
		${\begin{array}{c c}{\bf I}_2(p)\\p\geq5\end{array}}$& ${\begin{array}{c r}
			\Z_2  &p \text{ odd}\\	 
			\Z_2\oplus \Z_2 &p\text{ even}\\\end{array}}$ & ${\begin{array}{c r}
			0  &p \text{ odd}\\	 
			\Z_2 &p\text{ even}\\\end{array}}$ & ${\begin{array}{c r}
			\Z_2\oplus \Z_p  &p \text{ odd}\\	 
			\Z_2\oplus \Z_2 \oplus \Z_p&p\text{ even}\\\end{array}}$ \\	\hline&&&\\[-1em]
		${\bf F}_4$ & $\Z_2\oplus \Z_2$ &$\Z_2\oplus \Z_2$ &$\Z_2^5\oplus \Z_3^2 \oplus \Z_4$\\ \hline&&&\\[-1em]		
		${\bf H}_3$ & $\Z_2$ &$\Z_2$ & $\Z_2^3\oplus \Z_3 \oplus \Z_5$ \\ \hline&&&\\[-1em]
		${\bf H}_4$ & $\Z_2$ &$\Z_2$ &$\Z_2^2\oplus \Z_3 \oplus \Z_4 \oplus \Z_5$\\ \hline&&&\\[-1em]
		${\bf E}_6$ & $\Z_2$ &$\Z_2$ &$\Z_2^2\oplus \Z_3 \oplus \Z_4$\\ \hline&&&\\[-1em]
		${\bf E}_7$ & $\Z_2$ &$\Z_2$ &$\Z_2^2\oplus \Z_3 \oplus \Z_4$ \\\hline&&&\\[-1em]
		${\bf E}_8$ & $\Z_2$ &$\Z_2$ &$\Z_2^2\oplus \Z_3 \oplus \Z_4$ \\\hline
	\end{tabular}
	\caption{Homology of finite Coxeter groups.}\label{table:homology}
\end{table}

\bibliography{mybib}{}
\bibliographystyle{alpha}

\end{document}